\documentclass[pdflatex,sn-mathphys-num]{sn-jnl}
\usepackage{mathrsfs}
\usepackage{enumitem}
\usepackage{graphicx}%
\usepackage{multirow}%
\usepackage[all]{xy}
\usepackage{amsmath,amssymb,amsfonts}%
\usepackage{amsthm}%
\usepackage{mathrsfs}%
\usepackage[title]{appendix}%
\usepackage{xcolor}%
\usepackage{textcomp}%
\usepackage{manyfoot}%
\usepackage{booktabs}%
\usepackage{algorithm}%
\usepackage{algorithmicx}%
\usepackage{algpseudocode}%
\usepackage{listings}%
\usepackage{soul}
\usepackage{enumitem}
\usepackage{tikz-cd}
\theoremstyle{thmstyleone}%
\newtheorem{theorem}{Theorem}[section]%  meant for continuous numbers
\numberwithin{equation}{section}
\newtheorem{proposition}[theorem]{Proposition}% 

\theoremstyle{thmstyletwo}%
\newtheorem{example}[theorem]{Example}%
\newtheorem{remark}[theorem]{Remark}%

\theoremstyle{thmstylethree}%
\newtheorem{definition}[theorem]{Definition}%

\newtheorem{corollary}[theorem]{Corollary}
\begin{document}

    \title[Intermediate topological entropies for subsets of nonautonomous dynamical systems]{Intermediate topological entropies for subsets of nonautonomous dynamical systems}

%%=============================================================%%
%% GivenName	-> \fnm{Joergen W.}
%% Particle	-> \spfx{van der} -> surname prefix
%% FamilyName	-> \sur{Ploeg}
%% Suffix	-> \sfx{IV}
%% \author*[1,2]{\fnm{Joergen W.} \spfx{van der} \sur{Ploeg} 
%%  \sfx{IV}}\email{iauthor@gmail.com}
%%=============================================================%%

\author[1]{\fnm{Yujun} \sur{Ju}} \email{yjju@ctbu.edu.cn}

\affil[1]{\orgdiv{School of Mathematics and Statistics}, \orgname{Chongqing Technology and Business University}, \orgaddress{\city{Chongqing}, \postcode{400067}, \country{People's Republic of China}}}

%%==================================%%
%% Sample for unstructured abstract %%
%%==================================%%

\abstract{
Motivated by the notion of intermediate dimensions introduced by Falconer et al., we introduce a continuum of topological entropies that are intermediate between the (Bowen) topological entropy and the lower and upper capacity topological entropies.
This is achieved by restricting the families of allowable covers in the definition of topological entropy by requiring that the lengths of all strings used in a particular cover satisfy $ N \le n < N/\theta + 1 $, where $\theta \in [0,1] $ is a parameter.
When $\theta = 1 $, only covers using strings of the same length are allowed, and we recover the lower and upper capacity topological entropies; when $\theta = 0$, there are no restrictions, and the definition coincides with the topological entropy.
We first establish a quantitative inequality for the upper and lower intermediate topological entropies, which mirrors the corresponding result for intermediate dimensions. As a consequence, the intermediate topological entropies are continuous on $(0,1]$, though discontinuity may arise at $0$; an illustrative example is provided to demonstrate this phenomenon.
We then investigate several fundamental properties of the intermediate topological entropies for nonautonomous dynamical systems, including the power rule, monotonicity and product formulas. Finally, we derive an inequality relating intermediate entropies with respect to factor maps.}



\maketitle

\section{Introduction}
Topological entropy is one of the most fundamental invariants in dynamical systems, measuring the exponential growth rate of distinguishable orbits.
It was first introduced by Adler, Konheim and McAndrew \cite{adler1965topological} via open covers.
Bowen \cite{bowen1971entropy} and Dinaburg \cite{dinaburg1970correlation} later provided equivalent formulations using spanning and separated sets, and this definition mirrors the definition of the (upper) box dimension.
In 1973, Bowen \cite{bowen1973topological} defined the topological entropy on subsets resembling the definition of the Hausdorff dimension, which is called the Bowen topological entropy. Pesin \cite{pesin1997dimension} introduced a refined framework extending the classical Carath\'eodory construction, which is now commonly referred to as the Carath\'eodory-Pesin structure (C-P structure).
This framework characterizes the topological entropy and the lower and upper capacity topological entropies of noncompact sets, and also provides a unified approach to defining and studying quantities such as the Hausdorff and box dimensions, and the topological pressure in dynamical systems. For the packing dimension, its dynamical counterpart is the packing topological entropy, introduced by Feng and Huang \cite{feng2012variational}, who also established variational principles for both the Bowen and packing topological entropies by applying methods from geometric measure theory. Thus, the analogies between dimension and entropy are profound and far-reaching.

Falconer, Fraser and Kempton \cite{falconer2020intermediate} introduced a family of dimensions, called the intermediate dimensions, depending on a parameter $\theta \in [0,1]$.
These dimensions are defined by restricting the sizes of covering sets to lie within intervals of the form $[\delta^{1/\theta}, \delta]$.
The Hausdorff and box dimensions arise as the two extreme cases corresponding to $\theta = 0$ and $\theta = 1$, respectively.  Intermediate dimensions possess several useful properties, including continuity on $(0,1]$ (but not necessarily at $0$) and analogues of the mass distribution principle, Frostman's lemma and product formulas.
They also provide insight into the distribution of covering scales for sets whose Hausdorff and box dimensions differ, offering a refined description of geometric complexity. Moreover, since the intermediate dimensions are preserved under bi-Lipschitz mappings, they
provide another invariant for Lipschitz classification of sets. A very specific variant was used in \cite{kukavica2012estimate} to estimate the singular sets of partial differential equations.
Further developments and related results on intermediate dimensions can be found in \cite{Falconer2021Intermediate, banaji2022attainable, banaji2023intermediate,banaji2023generalised,burrell2022dimensions, daw2023fractal, falconer2022intermediate,feng2024intermediate}, among others.

Given the deep analogies between entropy and dimension theories, it is natural to seek dynamical counterparts of the intermediate dimensions, namely the \emph{intermediate topological entropies}.
The aim of this paper is to introduce and investigate the intermediate topological entropies for nonautonomous dynamical systems (NDSs for short).
The classical topological entropy for NDSs was first introduced and studied by Kolyada and Snoha \cite{kolyada1996topological}. Subsequently, Li \cite{li2015remarks} provided a concise overview of several notions of topological entropies for NDSs, including the classical, Bowen and Pesin topological entropies, and established the sufficient condition for their equivalence. Li and Ye \cite{li2023comparison} further discussed the interrelationships among the distance entropy, Bowen topological entropy, Pesin topological entropy, and the classical one, and established another equivalence condition for these quantities. Xu and Zhou \cite{xu2018variational}, as well as Zhang and Zhu \cite{zhang2023variational}, established variational principles for the Bowen and packing topological entropies of NDSs, respectively, by employing the approach developed by Feng and Huang. For more information on NDSs, the reader is referred to \cite{
kolyada1999topological,huang2008pressure,zhu2012entropy,canovas2013entropy,kuang2013fractal,kawan2013metric,kong2015topological,kawan2016some,Bis2018topological,xu2018variational,li2019topological,liu2020topological,rodrigues2022mean,shao2023topological,zhang2023variational,nazarian2024variational,wang2025entropies,li2025dimension,chen2025nonautonomous,chen2025nonautonomousII,yang2026bs} and references therein.

Motivated by the work of Falconer et al.~\cite{falconer2020intermediate} on intermediate dimensions, 
we define the lower and upper $\theta$-intermediate topological entropies for NDSs based on the idea of the C-P structure. 
This construction interpolates between the Bowen topological entropy and the classical topological entropy for NDSs by restricting the lengths of strings in admissible covers according to a parameter $\theta \in [0,1]$. 
When $\theta = 0$, both definitions reduce to the Bowen topological entropy, and when $\theta = 1$, the upper one coincides with the classical topological entropy. 
This definition provides a unified framework for describing dynamical complexity under varying covering scales and further deepens the connections between entropy and dimension theories. We note that the $\theta$-intermediate topological entropies are designed as dynamical counterparts of the intermediate dimensions, and it is expected that they share analogous structural properties. 
If such parallels did not appear, this would suggest that the construction is not properly aligned with its geometric motivation.

The paper is organized as follows. 
In Section~\ref{sec2}, we introduce the lower and upper $\theta$-intermediate topological entropies for NDSs. 
We then establish a quantitative inequality mirroring the corresponding result for intermediate dimensions, which leads to the continuity of the $\theta$-intermediate topological entropies on $(0,1]$, with possible discontinuities at $\theta=0$. 
Finally, by adapting existing arguments in the literature, we obtain Billingsley-type estimates for the $\theta$-intermediate topological entropies.
In Section~\ref{sec3}, we discuss several fundamental properties of the $\theta$-intermediate topological entropies, including closure stability, power rule, monotonicity and product formulas.
In Section~\ref{sec4}, the relations of the $\theta$-intermediate topological entropies between two topologically equisemiconjugate systems are studied, and we obtain an inequality for the $\theta$-intermediate entropies via factor maps.
In Section~\ref{ex}, we present an example of the $\theta$-intermediate topological entropy which is discontinuous at $\theta = 0$ but constant on $(0,1]$, thereby verifying the continuity inequality established in Section~\ref{sec2}.

\section{Intermediate topological entropies: definition and basic properties}\label{sec2}

In this section, we define and study the $\theta$-intermediate topological entropies for nonautonomous dynamical systems. 
After recalling the notions of intermediate dimensions and the classical topological entropy of NDSs, 
we establish several fundamental properties of intermediate topological entropies, as well as continuity with respect to~$\theta$, an equivalent formulation on compact metric spaces, 
and upper and lower estimates of intermediate topological entropies by local measure entropies.

\subsection{Review of intermediate dimensions}
We recall below the definition of the lower and upper \(\theta\)-intermediate dimensions, which was introduced and systematically studied in \cite{falconer2020intermediate}.
\begin{definition}
Let \( F \subseteq \mathbb{R}^{n} \) be bounded. For \( 0 \leq \theta \leq 1 \) the lower \( \theta \)-intermediate dimension of \( F \) is defined by
\[
\begin{aligned}
\underline{\dim}_{\theta} F = \inf \Bigl\{ & s \geq 0: \text{ for all } \varepsilon > 0 \text{ and all } \delta_0 > 0, \text{ there exists } 0 < \delta \leq \delta_0 \\
& \text{and a cover } \{ U_i \} \text{ of } F \text{ such that } \delta^{1/\theta} \leq |U_i| \leq \delta \text{ and } \sum |U_i|^s \leq \varepsilon \Bigr\}.
\end{aligned}
\]
Similarly, the upper \( \theta \)-intermediate dimension of \( F \) is defined by
\[
\begin{aligned}
\overline{\dim}_{\theta} F = \inf \Bigl\{& s \geq 0: \text{ for all } \varepsilon > 0 \text{ there exists } \delta_0 > 0 \text{ such that for all } 0 < \delta \leq \delta_0,\\
 &\text{there is a cover } \{ U_i \} \text{ of } F \text{ such that } \delta^{1/\theta} \leq |U_i| \leq \delta \text{ and } \sum |U_i|^s \leq \varepsilon \Bigr\}.
\end{aligned}
\]
\end{definition}
By convention, we set $\delta^{1/\theta}=0$ when $\theta=0$. It is immediate that
\[
\underline{\dim}_0 F = \overline{\dim}_0 F = \dim_{\mathrm{H}} F, \quad \underline{\dim}_1 F = \underline{\dim}_{\mathrm{B}} F \quad \text{and} \quad \overline{\dim}_1 F = \overline{\dim}_{\mathrm{B}} F,
\]
where $\dim_{\mathrm{H}} F$ denotes the Hausdorff dimension, 
and $\underline{\dim}_{\mathrm{B}} F$ and $\overline{\dim}_{\mathrm{B}} F$ 
represent the lower and upper box dimensions, respectively.
Furthermore, if $F\subseteq\mathbb{R}^n$ is bounded and $\theta\in[0,1]$, 
the $\theta$-intermediate dimensions satisfy the bounds between Hausdorff and box dimensions:
\[
0\le \dim_{\mathrm H}F\le \underline{\dim}_{\theta}F\le 
\overline{\dim}_{\theta}F\le\overline{\dim}_{\mathrm B}F\le n
~\text{and}~
0\le \underline{\dim}_{\theta}F\le \underline{\dim}_{\mathrm B}F\le n.
\]
As with box dimensions, when the lower and upper $\theta$-intermediate dimensions coincide we write simply 
$\dim_{\theta}F=\underline{\dim}_{\theta}F=\overline{\dim}_{\theta}F$ for the \(\theta\)-intermediate dimension of \(F\).

\subsection{Classical topological entropy of NDSs}
Following \cite{kolyada1996topological}, we recall the classical definition of topological entropy for nonautonomous dynamical systems and some of its equivalent formulations.
Let \( X \) be a compact topological space and let \( f_{1,\infty} = \{f_i\}_{i=1}^\infty \) be a sequence of continuous self-maps of \( X \). 
Denote by \(\mathbb{N}\) the set of positive integers.
For each \( i \in \mathbb{N} \), set \( f_i^0 = \text{id}_X \), the identity map of \( X \), and for each \( j \in \mathbb{N} \) define
\[
f_i^j = f_{i+(j - 1)} \circ \cdots \circ f_{i+1} \circ f_i, \quad f_i^{-j} = \left(f_i^j\right)^{-1} = f_i^{-1} \circ f_{i+1}^{-1} \circ \cdots \circ f_{i+(j - 1)}^{-1},
\]
where the above notation is applied to sets, and we do not assume that the maps $f_i$ are invertible.
Then we call \( (X, f_{1,\infty}) \) a nonautonomous dynamical system. If \( f_i = f \) for every \( i \geq 1 \), then \( (X, f_{1,\infty}) \) reduces to the autonomous dynamical system \( (X, f) \).

For open covers \( \mathcal{U}_1, \mathcal{U}_2, \ldots, \mathcal{U}_n \) of \( X \), we define their join as
\[
\bigvee_{i=1}^n \mathcal{U}_i = \mathcal{U}_1 \vee \mathcal{U}_2 \vee \cdots \vee \mathcal{U}_n = \left\{U_1 \cap U_2 \cap \cdots \cap U_n : U_i \in \mathcal{U}_i, i = 1, 2, \ldots, n\right\}.
\]
For a finite open cover \( \mathcal{U} \) of \( X \), we denote 
\[
f_i^{-j}(\mathcal{U}) = \left\{f_i^{-j}(U) : U \in \mathcal{U}\right\},\qquad
\mathcal{U}_i^n = \bigvee_{j=0}^{n - 1} f_i^{-j}(\mathcal{U}) 
\]
and by \( \mathcal{N}(\mathcal{U}) \) the minimal possible cardinality of a subcover chosen from \( \mathcal{U} \). 
Let $Z \subseteq X$ be a nonempty subset, not necessarily compact or invariant under $f_{1,\infty}$.
We denote by $\mathcal{U}|_Z$ the cover $\{U \cap Z : U \in \mathcal{U}\}$ of $Z$.

\begin{definition}
Let $X$ be a compact topological space and $f_{1,\infty}=\{f_i\}_{i=1}^\infty$ be a sequence of continuous self-maps on $X$.
Let $Z \subseteq X$ be a nonempty subset. The topological entropy of the sequence of maps \( f_{1,\infty} \) on the set \( Z \) is defined by
\[
h(f_{1,\infty}; Z) = \sup\left\{h(f_{1,\infty}; Z; \mathcal{U}) : \mathcal{U} \text{ is an open cover of } X\right\},
\]
where
\[
h(f_{1,\infty}; Z; \mathcal{U}) = \limsup_{n \to \infty} \frac{1}{n} \log \mathcal{N}\left(\left.\mathcal{U}_1^n\right|_Z\right).
\]
\end{definition}

The topological entropy can also be defined via separated sets and spanning sets.  
Assume \((X,d)\) is a compact metric space.  
For any \(U\subseteq X\) write \(|U|=\sup\{d(x,y):x,y\in U\}\) for its diameter, and for an open cover \(\mathcal U\) define \(|\mathcal U|=\max_{U\in\mathcal U}|U|\).

For each \(n\in\mathbb N\) set
\[
d_n(x,y)=\max_{0\le j\le n-1}d\bigl(f_1^j(x),f_1^j(y)\bigr),\qquad x,y\in X.
\]
Since \(X\) is compact, \(d_n\) is a metric equivalent to \(d\).  
Given $\varepsilon>0$ and $x\in X$, define the $(n,\varepsilon)$–Bowen ball of radius $\varepsilon$ centered at $x$ of length $n$ by
\[
B_n(x,\varepsilon)
   = \left\{y\in X : d_n(x,y)<\varepsilon \right\}.
\]

Let \(Z\subseteq X\) be nonempty.
A set \(E\subseteq X\) is an \((n,\varepsilon)\)-spanning set of \(Z\) if for every \(y\in Z\) there exists \(x\in E\) with \(d_n(x,y)\le\varepsilon\);
a set \(F\subseteq Z\) is an \((n,\varepsilon)\)-separated set of \(Z\) if \(x\neq y\) in \(F\) implies \(d_n(x,y)>\varepsilon\).
Write \(r_n(f_{1,\infty},Z,\varepsilon)\) for the minimal cardinality of an \((n,\varepsilon)\)-spanning set and \(s_n(f_{1,\infty},Z,\varepsilon)\) for the maximal cardinality of an \((n,\varepsilon)\)-separated set.

Then the topological entropy of \(f_{1,\infty}\) on \(Z\) can be expressed in terms of spanning and separated sets as
\[
h(f_{1,\infty};Z)
=\lim_{\varepsilon\to0}\,\limsup_{n\to\infty}\frac{1}{n}\log r_n(f_{1,\infty},Z,\varepsilon)
=\lim_{\varepsilon\to0}\,\limsup_{n\to\infty}\frac{1}{n}\log s_n(f_{1,\infty},Z,\varepsilon).
\]
In particular, if $f_i = f$ for all $i \ge 1$, then $h(f_{1,\infty};Z)$ coincides with the upper capacity topological entropy of $f$ restricted to $Z$, denoted by $h(f,Z)$.

\subsection{Intermediate topological entropies of NDSs}
Let $X$ be a compact topological space and $ f_{1,\infty}=\{ f_i\}_{i=1}^{\infty} $ be a sequence of continuous maps from $ X $ to itself.
For a finite open cover $\mathcal{U}$ of $X$ and \(m \in \mathbb{N}\), let
\[
\mathcal{S}_{m}(\mathcal{U}):=\left\{\textbf{U}=(U_{0},U_{1},\ldots,U_{m-1}):\textbf{U}\in \mathcal{U}^{m}\right\},
\]
where $\mathcal{U}^{m}=\prod\limits_{i=1}^{m}\mathcal{U}$. For any string $\textbf{U}\in \mathcal{S}_{m}(\mathcal{U})$, define the length of $\textbf{U}$ to be $m(\textbf{U}):=m$. 
If $k \in \mathbb{N}$ with $1 \le k \le m(\mathbf{U})$ and $0 \le a \le m(\mathbf{U}) - k$, we denote by
\[
\mathbf U|_{[a,a+k-1]} := (U_a,U_{a+1},\ldots,U_{a+k-1}) \in \mathcal S_k(\mathcal U)
\]
the substring of $\mathbf U$ of length $k$ starting at position $a$. 
In particular, $\mathbf U|_{[0,k-1]}$ is the initial truncation of length $k$. 

For a given string $\textbf{U}=(U_0,U_1,\ldots,U_{m-1}) \in \mathcal{S}_{m}(\mathcal{U})$, we associate the set
\[
 X_{f_{1,\infty}}(\textbf{U})=\left\{x\in X : f_1^{j}(x)\in U_j, j=0,1,\ldots,m-1\right\}.
\]
When no confusion arises we simply write \(X(\mathbf U)\) for \(X_{f_{1,\infty}}(\mathbf U)\),
and likewise omit the subscript \(f_{1,\infty}\) from the quantities \(M\), \(\underline m\), \(\overline m\) defined below.

For any subset $Z\subseteq X, \alpha \ge 0$ and $\theta\in[0,1]$, define
\[
\underline{m}(Z,\alpha,\mathcal{U},\theta)=\liminf_{N\to\infty}M(Z,\alpha,\mathcal{U},N,\theta),\quad 
\overline{m}(Z,\alpha,\mathcal{U},\theta)=\limsup_{N\to\infty}M(Z,\alpha,\mathcal{U},N,\theta),
\]
where
\[
M(Z,\alpha,\mathcal{U},N,\theta):=\inf_{\mathcal{G}}\left\{\sum\limits_{\textbf{U} \in \mathcal{G}}\exp(-\alpha m(\textbf{U}))\right\}
\]
and the infimum is taken over all finite or countable collections of strings $\mathcal{G} \subseteq \bigcup_{N \leq m < N/\theta+1} \mathcal{S}_{m}(\mathcal{U})$ such that \(\mathcal{G}\) covers \(Z\) $\left({\rm i.e.,} \bigcup_{\textbf{U}\in \mathcal{G}}X(\textbf{U})\supseteq Z\right)$.

It is easy to show that there are critical values of \(\underline{m}(Z,\alpha,\mathcal{U},\theta)\) and \(\overline{m}(Z,\alpha,\mathcal{U},\theta)\) as
\[\underline{h}_{\text{top}}(f_{1,\infty},Z,\mathcal{U},\theta)=\inf\left\{\alpha: \underline{m}(Z,\alpha,\mathcal{U},\theta)=0\right\}=\sup\{\alpha: \underline{m}(Z,\alpha,\mathcal{U},\theta)=\infty\},\]
\[\overline{h}_{\text{top}}(f_{1,\infty},Z,\mathcal{U},\theta)=\inf\{\alpha:\overline{m}(Z,\alpha,\mathcal{U},\theta)=0\}=\sup\{\alpha:\overline{m}(Z,\alpha,\mathcal{U},\theta)=\infty\}.\]

\begin{definition}
Let $X$ be a compact topological space and $f_{1,\infty}=\{f_i\}_{i=1}^\infty$ be a sequence of continuous self-maps on $X$. 
For $Z\subseteq X$, the lower and upper $\theta$-intermediate topological entropies of $f_{1,\infty}$ on $Z$ are defined by
\[
\underline{h}_{\mathrm{top}}(f_{1,\infty},Z,\theta)
=\sup_{\mathcal U}\,\underline{h}_{\mathrm{top}}(f_{1,\infty},Z,\mathcal U,\theta),\qquad
\overline{h}_{\mathrm{top}}(f_{1,\infty},Z,\theta)
=\sup_{\mathcal U}\,\overline{h}_{\mathrm{top}}(f_{1,\infty},Z,\mathcal U,\theta),
\]
where the supremum is over all finite open covers \(\mathcal U\) of $X$.  
If they coincide, we refer to the common value as the 
\(\theta\)-intermediate topological entropy of $f_{1,\infty}$ on $Z$, 
and denote it by \(h_{\mathrm{top}}(f_{1,\infty},Z,\theta)\).  
In the autonomous case $f_i=f$ for all $i \ge 1$, these reduce to 
\(\underline{h}_{\mathrm{top}}(f,Z,\theta)\), 
\(\overline{h}_{\mathrm{top}}(f,Z,\theta)\), 
and \(h_{\mathrm{top}}(f,Z,\theta)\), respectively.
\end{definition}

\begin{remark}
(i)  When \(\theta=0\), since \(M(Z,\alpha,\mathcal U,N,0)\) is non-decreasing with respect to \(N\), we have
\[
\underline{h}_{\mathrm{top}}(f_{1,\infty},Z,0)
=\overline{h}_{\mathrm{top}}(f_{1,\infty},Z,0).
\]
Following \cite{li2015remarks}, we call this the \emph{Pesin topological entropy},
which was originally introduced for compact metric spaces. 
In the same work, Li also introduced the Bowen topological entropy for NDSs, following Bowen's method \cite{bowen1973topological} for noncompact subsets of autonomous systems, and proved that the Pesin and Bowen topological entropies coincide. We therefore denote this common value by \(h_{\mathrm{top}}^{B}(f_{1,\infty},Z).\) 
For completeness, we briefly recall the definition of the Bowen topological entropy for NDSs.

Let \(\mathcal U\) be an open cover of \(X\) and \(K \subset X\). We write \(K \prec \mathcal U\) if \(K\) is contained in some element of \(\mathcal U\). Define
\[
n(f_{1,\infty}, \mathcal U, K)
:= \max\{n \in \mathbb{N} : f_1^i(K) \prec \mathcal U \text{ for all } 0 \le i < n\},
\]
with the convention that \(n(f_{1,\infty}, \mathcal U, K)=0\) if \(K \not\prec \mathcal U\), and \(n(f_{1,\infty}, \mathcal U, K)=+\infty\) if \(f_1^i(K) \prec \mathcal U\) for all \(i \in \mathbb{N}\).

Set
\[
B(f_{1,\infty}, \mathcal U, K):=\exp\big(-n(f_{1,\infty}, \mathcal U, K)\big),
\]
and for any countable collection \(\mathcal K=\{K_i\}_{i=1}^{\infty}\) and \(\lambda \in \mathbb{R}\), define
\[
B(f_{1,\infty}, \mathcal U, \mathcal K, \lambda)
:= \sum_{i=1}^{\infty} \big(B(f_{1,\infty}, \mathcal U, K_i)\big)^\lambda.
\]

For any given set \(Z \subset X\) (not necessarily compact or invariant), define
\[
m_{f_{1,\infty},\mathcal U}(\lambda,Z)
:= \lim_{\varepsilon\to0}
\inf\left\{
B(f_{1,\infty}, \mathcal U, \mathcal K, \lambda):
Z\subset\bigcup_{i=1}^{\infty} K_i,\ 
B(f_{1,\infty}, \mathcal U, K_i)<\varepsilon
\right\}.
\]

The Bowen topological entropy of \(f_{1,\infty}\) on \(Z\) is defined by
\begin{align*}
h_{\mathrm{top}}^{B}(f_{1,\infty},Z)
&:= \sup_{\mathcal U}
\inf\{\lambda: m_{f_{1,\infty},\mathcal U}(\lambda,Z)=0\} \\
&= \sup_{\mathcal U}
\sup\{\lambda: m_{f_{1,\infty},\mathcal U}(\lambda,Z)=+\infty\}.
\end{align*}

(ii) When \(\theta=1\), we have 
\[
M(Z,\alpha,\mathcal{U},N,1)=\inf\limits_{\mathcal{G}}\left\{\sum\limits_{\textbf{U} \in \mathcal{G}}\exp(-\alpha N)\right\}=\Lambda(Z,\mathcal{U},N)\exp(-\alpha N),
\]
where the infimum is taken over all finite or countable \(\mathcal G\subseteq\mathcal S_{N}(\mathcal U)\) with \(\bigcup_{\mathbf U\in\mathcal G}X(\mathbf U)\supseteq Z\), and
\(
\Lambda(Z,\mathcal U,N):=\inf_{\mathcal G}\operatorname{card}(\mathcal G).
\)
In this case one obtains
\[
\underline{h}_{\mathrm{top}}(f_{1,\infty},Z,\mathcal U,1)
=\liminf_{N\to\infty}\frac{1}{N}\log\Lambda(Z,\mathcal U,N),
\]
\[
\overline{h}_{\mathrm{top}}(f_{1,\infty},Z,\mathcal U,1)
=\limsup_{N\to\infty}\frac{1}{N}\log\Lambda(Z,\mathcal U,N),
\]
(see Pesin \cite[Theorem~2.2]{pesin1997dimension} for the proof in autonomous systems). 
Taking the supremum over all finite open covers \(\mathcal U\) of \(X\), we define the lower and upper capacity topological entropies of $ f_{1,\infty} $ on $ Z $ as
\[
\underline{Ch}(f_{1,\infty},Z):=\underline{h}_{\mathrm{top}}(f_{1,\infty},Z,1),\qquad
\overline{Ch}(f_{1,\infty},Z):=\overline{h}_{\mathrm{top}}(f_{1,\infty},Z,1).
\]
\end{remark}

\begin{corollary}
Let \( X \) be a compact topological space and \( f_{1,\infty} \) be a sequence of continuous self-maps of \( X \).
Then for any $Z\subseteq X$,
\[
\overline{Ch}(f_{1,\infty},Z)=h(f_{1,\infty};Z).
\]
\begin{proof}
Note that
\[
\Lambda(Z,\mathcal U,N)=\mathcal{N}\left(\left.\mathcal{U}_1^N\right|_Z\right).
\] 
Hence, by the definitions of \(\overline{Ch}(f_{1,\infty},Z)\) and \(h(f_{1,\infty};Z)\), the two quantities are equal.
\end{proof}
This result agrees with Corollary~1 in \cite{Bis2018topological} for compact metric spaces; 
see also the subsection below for the equivalent formulation in that setting.
\end{corollary}

Given two open covers \(\mathcal U\) and \(\mathcal V\) of \(X\), we say that \(\mathcal V\) is finer than \(\mathcal U\) (written \(\mathcal U\preceq\mathcal V\)) if every \(V\in\mathcal V\) is contained in some \(U\in\mathcal U\).
The following fundamental properties of lower and upper \(\theta\)–intermediate topological entropies can be verified directly from the definitions, so we omit the proofs for brevity.

\begin{proposition}\label{3.1}
Let $X$ be a compact topological space, and let $\mathcal{U}$ and $\mathcal{V}$ be open covers of $X$. 
For any $Z \subseteq X$ and $\theta \in [0,1]$,
\begin{itemize}
\item[(1)] $ \underline{h}_{\mathrm{top}}(f_{1,\infty},Z,\theta)\le \overline{h}_{\mathrm{top}}(f_{1,\infty},Z,\theta)$;
\item[(2)] If $Z_{1}\subseteq Z_{2} \subseteq X$, then 
\[
\overline{h}_{\mathrm{top}}(f_{1,\infty},Z_{1},\theta)\leq \overline{h}_{\mathrm{top}}(f_{1,\infty},Z_{2},\theta)\quad\text{and}\quad 
\underline{h}_{\mathrm{top}}(f_{1,\infty},Z_{1},\theta)\leq \underline{h}_{\mathrm{top}}(f_{1,\infty},Z_{2},\theta);
\]
\item[(3)] If $Z_{i}\subseteq X,\; i\geq1$ and $Z=\bigcup_{i\ge1}Z_{i}$, then 
\[
\overline{h}_{\mathrm{top}}(f_{1,\infty},Z,\theta)\ge
\sup_{i\ge1}\overline{h}_{\mathrm{top}}(f_{1,\infty},Z_{i},\theta)
\quad\text{and}\quad
\underline{h}_{\mathrm{top}}(f_{1,\infty},Z,\theta)\ge
\sup_{i\ge1}\underline{h}_{\mathrm{top}}(f_{1,\infty},Z_{i},\theta);
\]
\item[(4)] If $\mathcal{U}\preceq\mathcal{V}$, then 
\[
\overline{h}_{\mathrm{top}}(f_{1,\infty},Z,\mathcal{U},\theta)\leq 
\overline{h}_{\mathrm{top}}(f_{1,\infty},Z,\mathcal{V},\theta)
\quad\text{and}\quad
\underline{h}_{\mathrm{top}}(f_{1,\infty},Z,\mathcal{U},\theta)\leq 
\underline{h}_{\mathrm{top}}(f_{1,\infty},Z,\mathcal{V},\theta);
\]
\item[(5)] (\emph{Monotonicity in $\theta$})
If $0 \le \theta<\phi\le 1$, then
\[
\overline{h}_{\mathrm{top}}(f_{1,\infty},Z,\theta)\le 
\overline{h}_{\mathrm{top}}(f_{1,\infty},Z,\phi)
\quad\text{and}\quad
\underline{h}_{\mathrm{top}}(f_{1,\infty},Z,\theta)\le 
\underline{h}_{\mathrm{top}}(f_{1,\infty},Z,\phi);
\]

\item[(6)] (\emph{Finite stability})
Analogously to the upper \(\theta\)-intermediate dimension, the upper
\(\theta\)-intermediate topological entropy is finitely stable: for any 
\(Z_1,Z_2\subseteq X\),
\[
\overline{h}_{\mathrm{top}}(f_{1,\infty},Z_1 \cup Z_2,\theta)=\max\left\{
\overline{h}_{\mathrm{top}}(f_{1,\infty},Z_1,\theta),\overline{h}_{\mathrm{top}}(f_{1,\infty},Z_2,\theta)\right\}.
\]
\end{itemize}
\end{proposition}

\begin{remark}
As with the lower \(\theta\)-intermediate dimension, the lower \(\theta\)-intermediate topological entropy is in general not finitely stable. Moreover, neither the upper nor the lower \(\theta\)-intermediate topological entropy is countably stable.
\end{remark}

\subsection{Continuity with respect to $\theta$}
It is known that for a non-empty bounded subset \(F\subseteq\mathbb{R}^n\), the lower and upper
\(\theta\)-intermediate dimensions \(\underline{\dim}_{\theta}F\) and \(\overline{\dim}_{\theta}F\)
vary continuously for \(\theta\in(0,1]\), with possible discontinuities at \(\theta=0\)
(see Section~3.2 of \cite{falconer2020intermediate}). 
In particular, Proposition 14.2 of \cite{Falconer2021Intermediate} shows that for 
\(0<\theta<\phi\le1\) the upper \(\theta\)-intermediate dimensions satisfy
\[
\overline{\dim}_{\theta}F
\;\le\;
\overline{\dim}_{\phi}F
\;\le\;
\frac{\phi}{\theta}\,\overline{\dim}_{\theta}F,
\]
with analogous inequalities for the lower \(\theta\)-intermediate dimensions
(we note that \cite{Falconer2021Intermediate} also gives a second, slightly different 
inequality relating \(\theta\) and \(\phi\); here we only record the one relevant to our result). 
It is therefore natural to ask whether the lower and upper 
\(\theta\)-intermediate topological entropies 
\(\underline{h}_{\mathrm{top}}(f_{1,\infty},Z,\theta)\) and 
\(\overline{h}_{\mathrm{top}}(f_{1,\infty},Z,\theta)\) exhibit a similar behaviour. 
We show below that this is indeed the case: for each fixed \(Z\subseteq X\) they are continuous on 
\((0,1]\), with possible discontinuities at \(\theta=0\), exactly mirroring the situation for the 
intermediate dimensions. We provide a
simple example exhibiting discontinuity at \(\theta=0\) in Section \ref{ex}.

\begin{proposition}
Let \( X \) be a compact topological space and \( f_{1,\infty} \) be a sequence of continuous self-maps of \( X \).
Then for any nonempty $Z\subseteq X$ and \(0<\theta<\phi\le1\),
\[
\overline{h}_{\mathrm{top}}(f_{1,\infty},Z,\theta)
\;\le\;
\overline{h}_{\mathrm{top}}(f_{1,\infty},Z,\phi)
\;\le\;
\frac{\phi}{\theta}\,\overline{h}_{\mathrm{top}}(f_{1,\infty},Z,\theta),
\]
with the same inequalities for $\underline{h}_{\mathrm{top}}$.

\begin{proof}
The left-hand inequality follows from the monotonicity of \(\overline{h}_{\mathrm{top}}(f_{1,\infty},Z,\theta)\) in \(\theta\). To prove the right-hand inequality, fix \(0<\theta<\phi\le1\) and a finite open cover \(\mathcal U\) of \(X\).
Let \(s>\overline{h}_{\mathrm{top}}(f_{1,\infty},Z,\mathcal U,\theta)\) and \(\varepsilon>0\).
By definition, there exists \(N_0\) such that for all \(N\ge N_0\) we can choose a family \(\mathcal{G}\subseteq \bigcup_{N\le p < N/\theta+1}\mathcal{S}_p(\mathcal{U})\) covering \(Z\) with 
\[
\sum_{\mathbf{U}\in\mathcal{G}} \exp(-s\,m(\mathbf{U}))<\varepsilon.
\]
Split \(\mathcal{G}=\mathcal{G}_0 \cup \mathcal{G}_1\) where
\[
\mathcal{G}_0=\left\{\mathbf{U}\in\mathcal{G}: N\le m(\mathbf{U})<N/\phi+1\right\},\quad
\mathcal{G}_1=\left\{\mathbf{U}\in\mathcal{G}: N/\phi+1\le m(\mathbf{U})<N/\theta+1\right\}.
\]
Set \(q=\Big\lfloor \frac{N}{\phi}\Big\rfloor\), where \( \lfloor x \rfloor \) denotes the greatest integer less than
or equal to \(x\) and define
\[
\mathcal{G}_{1}^{*}:=\left\{\mathbf{U}^{(\phi)}=\mathbf{U}|_{[0,q-1]}:\mathbf{U}\in\mathcal{G}_{1}\right\}\subseteq
\mathcal{S}_q(\mathcal{U}).
\]
By the construction of the prefix $\mathbf{U}^{(\phi)}$, $X(\mathbf{U}^{(\phi)})\supseteq X(\mathbf{U})$,
so $\mathcal{G}_{1}^{*}$ together with $\mathcal{G}_{0}$ still covers $Z$.
Now put
\[
t_N=\frac{s\big(\frac{N}{\theta}+1\big)}{\big\lfloor \frac{N}{\phi}\big\rfloor}.
\]
Then for any \(\mathbf{U}\in\mathcal{G}_1\),
\[
t_{N}\Big\lfloor\frac{N}{\phi}\Big\rfloor = s\Big(\frac{N}{\theta}+1\Big) >s\,m(\mathbf{U}),
\]
and hence
\[
\exp\Bigl(-t_{N}\Bigl\lfloor \frac{N}{\phi}\Bigr\rfloor\Bigr)
<
\exp\bigl(-s\,m(\mathbf{U})\bigr).
\]
Therefore,
\begin{align*}
&\sum_{\mathbf{U} \in \mathcal{G}_0}\exp(-t_{N}\,m(\mathbf{U}))
   +\sum_{\mathbf{U}^{(\phi)} \in \mathcal{G}_1^*}\exp\Bigl(-t_{N}\,\Bigl\lfloor \tfrac{N}{\phi}\Bigr\rfloor\Bigr)\\
\leq & \sum_{\mathbf{U} \in \mathcal{G}_0}\exp(-s\, m(\mathbf{U}))
   +\sum_{\mathbf{U} \in \mathcal{G}_1}\exp(-s\, m(\mathbf{U})) \\
= & \sum_{\mathbf{U}\in\mathcal{G}}\exp(-s\, m(\mathbf{U})) < \varepsilon.
\end{align*}
It follows that for all \(N\ge N_0\),
\[
M(Z,t_N,\mathcal{U},N,\phi) < \varepsilon.
\]
Since \(t_N\to \frac{\phi}{\theta}\,s\) as \(N\to\infty\), fix any \(t>\frac{\phi}{\theta}\,s\); then
there exists \(N_1\) such that \(t>t_N\) for all \(N\ge N_1\).
Hence for all \(N\ge \max\{N_0,N_1\}\),
\[
M(Z,t,\mathcal{U},N,\phi) \leq M(Z,t_N,\mathcal{U},N,\phi) < \varepsilon.
\]
Taking the upper limit as \(N\to\infty\) gives
\(\overline m(Z,t,\mathcal U,\phi)=0\),
and consequently
\[
\overline{h}_{\mathrm{top}}(f_{1,\infty},Z,\mathcal{U},\phi)\ \le\ t.
\]
Since \(t>\frac{\phi}{\theta}\,s\) is arbitrary, we conclude
\[
\overline{h}_{\mathrm{top}}(f_{1,\infty},Z,\mathcal{U},\phi)\ \le\ \frac{\phi}{\theta}\,s.
\]
It follows that
\[
\overline{h}_{\mathrm{top}}(f_{1,\infty},Z,\mathcal{U},\phi)\ \le\ \frac{\phi}{\theta}\,
\overline{h}_{\mathrm{top}}(f_{1,\infty},Z,\mathcal{U},\theta).
\]
Therefore
\[
\overline{h}_{\mathrm{top}}(f_{1,\infty},Z,\phi)\ \le\ \frac{\phi}{\theta}\,
\overline{h}_{\mathrm{top}}(f_{1,\infty},Z,\theta).
\]
The corresponding statement for the lower \(\theta\)–intermediate topological entropy follows by exactly the same argument, taking covers of \(Z\) by strings with \(N\le m(\mathbf U)<N/\theta+1\) for arbitrarily large \(N\).
\end{proof}
\end{proposition}

\begin{corollary}
The maps \(\theta\mapsto \underline{h}_{\mathrm{top}}(f_{1,\infty},Z,\theta)\) and
\(\theta\mapsto \overline{h}_{\mathrm{top}}(f_{1,\infty},Z,\theta)\) are continuous for \(\theta\in(0,1]\).
\end{corollary}

\subsection{Equivalent definition on compact metric spaces}

We now present an equivalent definition of the lower and upper \(\theta\)-intermediate topological entropies for subsets of a compact metric space $(X,d)$. 
For any $\alpha \geq 0$, $N \in \mathbb{N}$, $\delta>0$ and $\theta\in [0,1]$, define
\[
\begin{aligned}
M(Z, \alpha, \delta,N, \theta) 
= \inf \Biggl\{
  &\sum_i \exp\bigl(-\alpha\,n_i\bigr): \\[0.7em]
  &\Bigl.\;
   \bigcup_i B_{n_i}(x_i,\delta)\supseteq Z,~
   x_i\in X,~N \le n_i  < N/\theta+1 \text{ for all }i
\Biggr\}.
\end{aligned}
\]

Let
\[\underline{m}(Z, \alpha, \delta,\theta)=\liminf_{N \rightarrow \infty} M(Z, \alpha, \delta, N,\theta),\]
\[\overline{m}(Z, \alpha, \delta,\theta)=\limsup_{N \rightarrow \infty} M(Z, \alpha, \delta, N,\theta).\]

We define the lower and upper \(\theta\)-intermediate topological entropies of \(Z\) relative to \(\delta\) by
\[\underline{h}_{\mathrm{top}}(f_{1,\infty},Z,\delta,\theta)=\inf \left\{\alpha: \underline{m}(Z, \alpha, \delta,\theta)=0\right\}=\sup \left\{\alpha: \underline{m}(Z, \alpha, \delta,\theta)=\infty\right\},
\]
\[\overline{h}_{\mathrm{top}}(f_{1,\infty},Z,\delta,\theta)=\inf \{\alpha: \overline{m}(Z, \alpha, \delta,\theta)=0\}=\sup \{\alpha: \overline{m}(Z, \alpha, \delta,\theta)=\infty\}.\]

It is easy to see that 
\(\underline{h}_{\mathrm{top}}(f_{1,\infty},Z,\delta,\theta)\)
is nondecreasing as \(\delta\) decreases (and so is \(\overline{h}_{\mathrm{top}}\)).

\begin{theorem}\label{thm:limit-entropy}
Let $(X,d)$ be a compact metric space and $f_{1,\infty}$ be a sequence of continuous self-maps of $X$. 
For any set $Z \subseteq X$ and $\theta \in [0,1]$, the following limits exist:
\[
\underline{h}_{\mathrm{top}}(f_{1,\infty},Z,\theta)
=\lim_{|\mathcal{U}|\to 0}\underline{h}_{\mathrm{top}}(f_{1,\infty},Z,\mathcal{U},\theta),
~
\overline{h}_{\mathrm{top}}(f_{1,\infty},Z,\theta)
=\lim_{|\mathcal{U}|\to 0}\overline{h}_{\mathrm{top}}(f_{1,\infty},Z,\mathcal{U},\theta).
\]
\end{theorem}
\begin{proof}
We give the proof for $\underline{h}_{\mathrm{top}}$; the argument for $\overline{h}_{\mathrm{top}}$ is identical.
By definition,
\[
\underline{h}_{\mathrm{top}}(f_{1,\infty},Z,\theta)
=\sup_{\mathcal{U}}\underline{h}_{\mathrm{top}}(f_{1,\infty},Z,\mathcal{U},\theta),
\]
where the supremum runs over all finite open covers $\mathcal{U}$ of $X$.
Thus
\[
\liminf_{|\mathcal{U}|\to 0}\underline{h}_{\mathrm{top}}(f_{1,\infty},Z,\mathcal{U},\theta)
\le
\limsup_{|\mathcal{U}|\to 0}\underline{h}_{\mathrm{top}}(f_{1,\infty},Z,\mathcal{U},\theta)
\le
\underline{h}_{\mathrm{top}}(f_{1,\infty},Z,\theta).
\]

Conversely, let $\mathcal{V}$ be a finite open cover of $X$ with Lebesgue number $\delta>0$.
If $\mathcal{U}$ is a finite open cover with $|\mathcal{U}|<\delta$, then $\mathcal{V}\preceq\mathcal{U}$.
By monotonicity with respect to refinement,
\[
\underline{h}_{\mathrm{top}}(f_{1,\infty},Z,\mathcal{V},\theta)
\le
\underline{h}_{\mathrm{top}}(f_{1,\infty},Z,\mathcal{U},\theta).
\]
Taking $\liminf_{|\mathcal{U}|\to 0}$ gives
\[
\underline{h}_{\mathrm{top}}(f_{1,\infty},Z,\mathcal{V},\theta)
\le
\liminf_{|\mathcal{U}|\to 0}\underline{h}_{\mathrm{top}}(f_{1,\infty},Z,\mathcal{U},\theta).
\]
Taking the supremum over all $\mathcal{V}$ yields
\[
\underline{h}_{\mathrm{top}}(f_{1,\infty},Z,\theta)
\le
\liminf_{|\mathcal{U}|\to 0}\underline{h}_{\mathrm{top}}(f_{1,\infty},Z,\mathcal{U},\theta).
\]
Combining with the first inequality we obtain
\[
\lim_{|\mathcal{U}|\to 0}\underline{h}_{\mathrm{top}}(f_{1,\infty},Z,\mathcal{U},\theta)
=
\underline{h}_{\mathrm{top}}(f_{1,\infty},Z,\theta).
\]
\end{proof}

\begin{theorem}
Let $(X,d)$ be a compact metric space and $f_{1,\infty}$ be a sequence of continuous self-maps of $X$. 
For any set $Z \subseteq X$ and $\theta \in [0,1]$, the following limits exist:
\[
\underline{h}_{\mathrm{top}}(f_{1,\infty},Z,\theta)
=\lim_{\delta \to 0}\underline{h}_{\mathrm{top}}(f_{1,\infty},Z,\delta,\theta),
\qquad
\overline{h}_{\mathrm{top}}(f_{1,\infty},Z,\theta)
=\lim_{\delta \to 0}\overline{h}_{\mathrm{top}}(f_{1,\infty},Z,\delta,\theta).
\]
\end{theorem}
\begin{proof}
Let \(\mathcal{U}\) be a finite open cover of \(X\) and \(\delta(\mathcal{U})\) be its Lebesgue number. 
It is easy to see that for every \(x \in X,\) if \(x \in X(\mathbf{U})\) for some 
\(\mathbf{U} \in \mathcal{S}_{n}(\mathcal{U})\), then
\[
X(\mathbf{U}) \subseteq B_{n}(x,2|\mathcal{U}|).
\]
Moreover, for each \(x \in X\), there exists \(\mathbf{U}_x \in \mathcal{S}_n(\mathcal{U})\) such that
\[
B_{n}\!\left(x,\frac{1}{2}\delta(\mathcal{U})\right) \subseteq X(\mathbf{U}_x).
\]
Thus
\[
M(Z,\alpha,2|\mathcal{U}|,N,\theta)
\le 
M(Z,\alpha,\mathcal{U},N,\theta)
\le 
M\!\left(Z,\alpha,\frac{1}{2}\delta(\mathcal{U}),N,\theta\right).
\]
This implies
\[
\underline{h}_{\mathrm{top}}\left(f_{1,\infty},Z,2|\mathcal{U}|,\theta\right)
\le
\underline{h}_{\mathrm{top}}\left(f_{1,\infty},Z,\mathcal{U},\theta\right)
\le
\underline{h}_{\mathrm{top}}\left(f_{1,\infty},Z,\frac{1}{2}\delta(\mathcal{U}),\theta\right),
\]
and similarly for \(\overline{h}_{\mathrm{top}}\). Letting \(|\mathcal{U}|\to0\) (hence \(\delta(\mathcal{U})\to0\))
yields the desired limits.
\end{proof}

\subsection{Billingsley-type estimates for $\theta$-intermediate topological entropies}

 Ma and Wen showed that Bowen topological entropy can be determined via local measure entropies in \cite[Theorem~1]{ma2008billingsley},  which is an analogue of Billingsley's theorem for the Hausdorff dimension. Following their work, Ju and Yang~\cite{ju2021pesin} extended this result to nonautonomous dynamical systems, while Bi\'s~\cite{Bis2018topological} independently obtained a similar relation for the classical topological entropy \(h(f_{1,\infty};Z)\).
Both results provide upper and lower estimates of topological entropies via local measure entropies.
In this section, we establish an analogous Billingsley-type result for the \(\theta\)-intermediate topological entropies of NDSs.

Let $(X,f_{1,\infty})$ be an NDS on a compact metric space $(X,d)$,
and let $\mathcal{M}(X)$ denote the set of Borel probability measures on $X$. A metric space is called boundedly compact if all bounded closed subsets are compact.
We say that a Borel probability measure \(\mu\) on a compact metric space \((X, d)\) is strictly positive if \(\mu(U) > 0\) for any nonempty open set \(U \subseteq X\). To establish our result, we first recall the relevant definitions and auxiliary theorems concerning local measure entropies and their connections with topological entropies. 

\begin{definition}
Following Bi\'s~\cite{Bis2018topological}, for \(\mu \in \mathcal{M}(X)\) and \(x \in X\), the quantity
\[
\overline{h}_{\mu, f_{1,\infty}}(x) = \lim_{\epsilon \to 0} \limsup_{n \to \infty} -\frac{1}{n} \log \mu(B_n(x, \epsilon))
\]
is called a local upper \( \mu \)-measure entropy at the point \( x \), with respect to \( f_{1,\infty} = \{f_n\}_{n=1}^\infty \), while the quantity
\[
\underline{h}_{\mu, f_{1,\infty}}(x) = \lim_{\epsilon \to 0} \liminf_{n \to \infty} -\frac{1}{n} \log \mu(B_n(x, \epsilon))
\]
is called a local lower \( \mu \)-measure entropy at the point \( x \), with respect to \( f_{1,\infty} = \{f_n\}_{n=1}^\infty \).
\end{definition}

The following theorem is due to Ju and Yang~\cite[Theorem~3.2]{ju2021pesin}. It provides a lower bound for the Bowen topological entropy in terms of local lower measure entropies.

\begin{theorem}
Let \((X,d)\) be a compact metric space and \( f_{1,\infty} \) be a sequence of continuous self-maps of \( X \). 
Then for any Borel subset $Z\subseteq X$, $s\in(0,\infty)$ and 
$\mu\in\mathcal{M}(X)$, if
\[
\underline{h}_{\mu,f_{1,\infty}}(x)\ge s
\quad \text{for all }x\in Z
~\text{and}~ \mu(Z)>0 \quad
\text{then} \quad h_{\mathrm{top}}^{B}(f_{1,\infty},Z)\ge s.
\]
\end{theorem}

Similarly, Bi\'s~\cite[Theorem~2]{Bis2018topological} obtained the following upper estimate for the upper capacity topological entropy by local upper measure entropies.

\begin{theorem}
Let $(X,d)$ be a boundedly compact metric space that admits a strictly positive Borel probability measure $\mu$, and let $f_{1,\infty}=\{f_n\}_{n\ge1}$ be a sequence of continuous self-maps of $X$.
Then for any Borel subset $Z\subseteq X$ and $s\in(0,\infty)$,  
if
\[
\overline{h}_{\mu,f_{1,\infty}}(x)\le s
\quad \text{for all }x\in Z \quad
\text{then} \quad h_{\mathrm{top}}(f_{1,\infty};Z)\le s.
\]
\end{theorem}

By the monotonicity of the $\theta$-intermediate topological entropy with respect to~$\theta$
and the two preceding theorems, we can immediately obtain the following result.

\begin{theorem}
Let \((X,d)\) be a compact metric space and \( f_{1,\infty} \) be a sequence of continuous self-maps of \( X \). 
Then for any Borel subset $Z\subseteq X$, $s\in(0,\infty)$ and $\theta\in[0,1]$, the following hold:
\begin{enumerate}
\item[(1)] If the space $X$ admits a strictly positive Borel probability
measure $\mu$ such that
\[
\overline{h}_{\mu,f_{1,\infty}}(x)\le s
\quad \text{for all }x\in Z \quad
\text{then} \quad \overline{h}_{\mathrm{top}}(f_{1,\infty},Z,\theta)\le s.
\]

\item[(2)] For any $\mu\in\mathcal{M}(X)$, if
\[
\underline{h}_{\mu,f_{1,\infty}}(x)\ge s
\quad \text{for all }x\in Z
~\text{and } \mu(Z)>0 \quad
\text{then} \quad \underline{h}_{\mathrm{top}}(f_{1,\infty},Z,\theta)\ge s.
\]
\end{enumerate}
\end{theorem}

\section{Further properties of intermediate topological entropies}\label{sec3}

In this section, we establish several fundamental properties of the intermediate topological entropies for $(X, f_{1,\infty})$, which make their computation more accessible.
Many of the results obtained in this section can be viewed as extensions of known properties of classical topological entropy \cite{kolyada1996topological} and Pesin topological entropy \cite{ju2021pesin} for nonautonomous dynamical systems.

\begin{proposition}\label{sc}
Let \((X,d)\) be a compact metric space and \( f_{1,\infty} \) be a sequence of continuous self-maps of \( X \). 
Then, for any \(Z\subseteq X\) and \(\theta\in(0,1]\),
\[
\underline{h}_{\mathrm{top}}(f_{1,\infty},\overline{Z},\theta)
=\underline{h}_{\mathrm{top}}(f_{1,\infty},Z,\theta),
\qquad
\overline{h}_{\mathrm{top}}(f_{1,\infty},\overline{Z},\theta)
=\overline{h}_{\mathrm{top}}(f_{1,\infty},Z,\theta).
\]
\begin{proof}
Since \(Z\subseteq \overline{Z}\), by Proposition \ref{3.1}(2), it suffices to prove the following inequalities 
\[ 
\underline h_{\mathrm{top}}(f_{1,\infty},\overline Z,\theta) \le \underline h_{\mathrm{top}}(f_{1,\infty},Z,\theta),\quad \overline h_{\mathrm{top}}(f_{1,\infty},\overline Z,\theta) \le \overline h_{\mathrm{top}}(f_{1,\infty},Z,\theta). 
\]

Fix \(N\in\mathbb N\) and note that for each integer \(0\le j< N/\theta+1<\infty\), 
the map \(f_1^j\) is uniformly continuous on the compact space \(X\).
Hence for every \(\delta>0\) there exists \(\eta_j>0\) such that 
\[
d(x,y)<\eta_j \Rightarrow d\bigl(f_1^j(x),f_1^j(y)\bigr)<\tfrac{\delta}{2}.
\]

Set \(\varepsilon_N = \min_{0 \leq j < N/\theta + 1} \eta_j > 0\). Then for every \(n\in [N,N/\theta+1)\), \(d(x,y)<\varepsilon_N\) implies \(d_n(x,y)<\delta/2\).

Now, let \(\left\{B_{n_i}(x_i, \delta/2)\right\}_i\) be any cover of \(Z\) with \(N \leq n_i < N/\theta + 1\). For any \(y \in \overline{Z}\), select \(z \in Z\) such that \(d(y,z) < \varepsilon_N\), and choose \(i\) with \(z \in B_{n_i}(x_i, \delta/2)\). Then,
\[
d_{n_i}(y, x_i) \leq d_{n_i}(y, z) + d_{n_i}(z, x_i) < \frac{\delta}{2} + \frac{\delta}{2} = \delta,
\]
so \(y \in B_{n_i}(x_i, \delta)\). Hence, \(\{B_{n_i}(x_i, \delta)\}_i\) covers \(\overline{Z}\) and consequently,
\[
M(\overline Z,\alpha,\delta,N,\theta) \;\le\; M\!\left(Z,\alpha,\tfrac{\delta}{2},N,\theta\right).
\]
Taking the \(\liminf\) or \(\limsup\) as \(N \to \infty\) gives
\[
\underline m(\overline Z,\alpha,\delta,\theta)
\le
\underline m\left(Z,\alpha,\tfrac{\delta}{2},\theta \right),
\qquad
\overline m(\overline Z,\alpha,\delta,\theta)
\le
\overline m\left(Z,\alpha,\tfrac{\delta}{2},\theta \right).
\]
Thus,
\[
\underline h_{\mathrm{top}}(f_{1,\infty},\overline Z,\delta,\theta)
\le
\underline h_{\mathrm{top}}\left(f_{1,\infty},Z,\tfrac{\delta}{2},\theta \right),
\qquad
\overline h_{\mathrm{top}}\left(f_{1,\infty},\overline Z,\delta,\theta \right)
\le
\overline h_{\mathrm{top}}\left(f_{1,\infty},Z,\tfrac{\delta}{2},\theta \right).
\]
Finally, letting \(\delta \to 0\) completes the proof.
\end{proof}
\end{proposition}

For $\theta\in(0,1]$, this result shows that the $\theta$-intermediate topological entropies 
are stable under closure, so one may restrict attention to compact (closed) subsets in computations.  
However, this closure stability is incompatible with the countable stability that occurs at $\theta=0$, 
indicating a fundamental difference between the two cases.

\begin{proposition}\label{pr-theta}
Let \( X \) be a compact topological space and \( f_{1,\infty} \) be a sequence of continuous self-maps of \( X \).
Then for any $Z\subseteq X, \theta \in [0,1]$ and $m\in\mathbb N$,
\[
\underline{h}_{\mathrm{top}}\left(f_{1,\infty}^m,Z,\theta \right)\le m\,\underline{h}_{\mathrm{top}}\left(f_{1,\infty},Z,\theta \right)
\quad\text{and}\quad
\overline{h}_{\mathrm{top}}\left(f_{1,\infty}^m,Z,\theta \right)\le m\,\overline{h}_{\mathrm{top}}\left(f_{1,\infty},Z,\theta \right),
\]
where $f_{1,\infty}^m=\left\{f_{im+1}^{m}\right\}_{i=0}^\infty$.
\end{proposition}
\begin{proof}
Since the case $\theta=0$ was proved in \cite[Theorem~3.3]{ju2021pesin} for compact metric spaces and the proof remains valid here, we consider only $\theta\in(0,1]$ henceforth.

Fix a finite open cover \(\mathcal U\) of \(X\) and let \(\alpha\ge0\).
Write \(k=mN+r\) with \(N\in\mathbb N\) and \(r\in\{0,\ldots,m-1\}\).
The proof proceeds in three steps.

\smallskip
\emph{Step~1.}
Consider any family
\[
\mathcal G_{k,f_{1,\infty}}
   \subseteq\bigcup_{k\le p<k/\theta+1}\mathcal S_p(\mathcal U)
\]
that covers \(Z\) with respect to \(f_{1,\infty}\).
For each \(\mathbf U=(U_0,\ldots,U_{p-1})\in\mathcal G_{k,f_{1,\infty}}\),
define \(\mathbf U^*:=\mathbf U|_{[0,p-r-1]}\).
Then \(m(\mathbf U^*) \in [mN,(mN+r)/\theta+1-r)\) and
\(X_{f_{1,\infty}}(\mathbf U)\subseteq X_{f_{1,\infty}}(\mathbf U^*)\).
Define
\[
\mathcal G_{k,f_{1,\infty}}^{*}
   :=\left\{\mathbf U^*:\mathbf U\in\mathcal G_{k,f_{1,\infty}}\right\}.
\]
Hence
\begin{equation}\label{eq:s1-en}
\sum_{\mathbf U\in\mathcal G_{k,f_{1,\infty}}}e^{-\alpha m(\mathbf U)}
   \ge e^{-\alpha r}
      \sum_{\mathbf U^* \in \mathcal G_{k,f_{1,\infty}}^{*}}e^{-\alpha m(\mathbf U^*)}.
\end{equation}

\emph{Step~2.}
If \(m(\mathbf U^*)<mN/\theta+1\), set \(\mathbf W(\mathbf U^*):=\mathbf U^*\);
otherwise, truncate further by setting
\[
q:=\Big\lfloor\frac{mN}{\theta}\Big\rfloor,
\qquad
\mathbf W(\mathbf U^*):=\mathbf U^*|_{[0,q-1]}.
\]
Then \(m\left(\mathbf W(\mathbf U^*)\right)\in[mN,mN/\theta+1)\) and
\(X_{f_{1,\infty}}(\mathbf U^*)\subseteq X_{f_{1,\infty}}\left(\mathbf W (\mathbf U^*)\right)\).
When the additional truncation occurs, we have 
\[
0
\le m(\mathbf U^*)-m\left(\mathbf W (\mathbf U^*)\right)
< \frac{mN+r}{\theta}+1-r-\Bigl\lfloor\frac{mN}{\theta}\Bigr\rfloor
< \frac{r}{\theta}+2-r <\frac{m}{\theta}+2-r.
\]
Define
\[
\mathcal G_{mN,f_{1,\infty}}
   :=\left\{\mathbf W(\mathbf U^*):\mathbf U^* \in \mathcal G_{k,f_{1,\infty}}^{*} \right\}.
\]
Then
\begin{equation}\label{eq:s2-en}
\sum_{\mathbf U^*\in\mathcal G_{k,f_{1,\infty}}^{*}}e^{-\alpha m(\mathbf U^*)}
   \ge e^{-\alpha\left(\frac{m}{\theta}+2-r\right)}
      \sum_{\mathbf W\in\mathcal G_{mN,f_{1,\infty}}}e^{-\alpha m(\mathbf W)}.
\end{equation}

\emph{Step~3.}
For each \(\mathbf W=(W_0,\ldots,W_{t-1})\in\mathcal G_{mN,f_{1,\infty}}\),
set \(n:=\lfloor t/m\rfloor\) and define
\[
\mathbf V(\mathbf W):=(V_0,\ldots,V_{n-1}),\qquad V_j:=W_{jm}.
\]
Then \(m\left(\mathbf V(\mathbf W)\right) \in [N,N/\theta+1)\) and
\(X_{f_{1,\infty}}(\mathbf W)\subseteq X_{f_{1,\infty}^m}\left(\mathbf V(\mathbf W)\right)\),
so that
\[
\mathcal G_{N,f_{1,\infty}^m}
   :=\left\{\mathbf V(\mathbf W):\mathbf W\in\mathcal G_{mN,f_{1,\infty}}\right\}
\]
is still a cover of \(Z\).
Since \(n=\lfloor t/m \rfloor\), we have \(mn \le t < m(n+1)\).
Hence
\[
e^{-\alpha t} > e^{-\alpha(mn+m)} = e^{-\alpha m}\,e^{-\alpha m n}.
\]
Summing over all \(\mathbf W\) yields
\begin{equation}\label{eq:s3-en}
\sum_{\mathbf W\in\mathcal G_{mN,f_{1,\infty}}}e^{-\alpha m(\mathbf W)}
   \ge e^{-\alpha m}
      \sum_{\mathbf V\in\mathcal G_{N,f_{1,\infty}^m}}
      e^{-\alpha m\,m(\mathbf V)}.
\end{equation}

Combining \eqref{eq:s1-en}–\eqref{eq:s3-en}, we obtain
\[
\sum_{\mathbf U\in\mathcal G_{k,f_{1,\infty}}}e^{-\alpha m(\mathbf U)}
   \ge e^{-\alpha (\frac{m}{\theta}+2+m)}
      \sum_{\mathbf V\in\mathcal G_{N,f_{1,\infty}^m}}
      e^{-\alpha m\,m(\mathbf V)},
\]
which implies
\[
M_{f_{1,\infty}}(Z,\alpha,\mathcal U,k,\theta)
   \ge 
   e^{-\alpha (\frac{m}{\theta}+2+m)}
   M_{f_{1,\infty}^m}(Z,\alpha m,\mathcal U,N,\theta).
\]

Taking $\liminf_{N\to\infty}$ (resp. $\limsup_{N\to\infty}$) with $k=mN+r$, we obtain
\[
\underline{m}_{f_{1,\infty}}(Z,\alpha,\mathcal U,\theta)
   \ge e^{-\alpha (\frac{m}{\theta}+2+m)}
      \underline{m}_{f_{1,\infty}^m}(Z,\alpha m,\mathcal U,\theta),
\]
\[
\overline{m}_{f_{1,\infty}}(Z,\alpha,\mathcal U,\theta)
   \ge e^{-\alpha (\frac{m}{\theta}+2+m)}
      \overline{m}_{f_{1,\infty}^m}(Z,\alpha m,\mathcal U,\theta).
\]

It follows that
\[
\underline{h}_{\mathrm{top}}(f_{1,\infty}^m,Z,\mathcal U,\theta)
\le m\,\underline{h}_{\mathrm{top}}(f_{1,\infty},Z,\mathcal U,\theta),
\quad
\overline{h}_{\mathrm{top}}(f_{1,\infty}^m,Z,\mathcal U,\theta)
\le m\,\overline{h}_{\mathrm{top}}(f_{1,\infty},Z,\mathcal U,\theta).
\]
Finally, taking the supremum over all finite open covers yields the desired inequalities.
\end{proof}

\begin{proposition}\label{pr1-theta}
Let \( X \) be a compact topological space and \( f_{1,\infty} \) be a sequence of continuous self-maps of \( X \).
Assume $f_{1,\infty}$ is periodic with period $m\in\mathbb N$, i.e.\ $f_{im+j}=f_j$ for all $i\ge0$ and $1\le j\le m$.
Then for any $Z\subseteq X \text{ and } \theta \in [0,1]$,
\[
\underline{h}_{\mathrm{top}}(f_{1,\infty}^m,Z,\theta)=
m\,\underline{h}_{\mathrm{top}}(f_{1,\infty},Z,\theta),
\qquad
\overline{h}_{\mathrm{top}}(f_{1,\infty}^m,Z,\theta)=
m\,\overline{h}_{\mathrm{top}}(f_{1,\infty},Z,\theta).
\]
\end{proposition}
\begin{proof}
The inequalities  
\[
\underline{h}_{\mathrm{top}}(f_{1,\infty}^m,Z,\theta)\le
m\,\underline{h}_{\mathrm{top}}(f_{1,\infty},Z,\theta)
\quad\text{and}\quad
\overline{h}_{\mathrm{top}}(f_{1,\infty}^m,Z,\theta)\le
m\,\overline{h}_{\mathrm{top}}(f_{1,\infty},Z,\theta)
\]
follow from Proposition~\ref{pr-theta}.  
Hence we only need to prove the reverse inequalities.
Fix a finite open cover \(\mathcal U\) of \(X\) and define
\[
\mathcal V := \mathcal U \vee f_1^{-1}\mathcal U \vee \cdots \vee f_1^{-(m-1)}\mathcal U.
\]
The argument proceeds in three steps.

\emph{Step~1.}
For any \(N\in\mathbb N\), write \(k=mN+r\) with \(r\in\{0,\ldots,m-1\}\) and let 
\[
\mathcal G_{N,f_{1,\infty}^m} \subseteq \bigcup_{N\le n<N/\theta+1}\mathcal S_n(\mathcal V)
\]
be a cover of \(Z\) with respect to \(f_{1,\infty}^m\).
For each \(\mathbf V=(V_0,\ldots,V_{n-1})\in\mathcal G_{N,f_{1,\infty}^m}\), write
\[
V_j=U_{jm}\cap f_1^{-1}U_{jm+1}\cap\cdots\cap f_1^{-(m-1)}U_{jm+m-1}\quad(U_{jm+\ell}\in\mathcal U),
\]
and concatenate to obtain \(\mathbf U(\mathbf V)=(U_0,\ldots,U_{p-1})\) with \(p=mn\).
By the periodicity \(f_{im+j}=f_j\), we have
\[
X_{f_{1,\infty}^m}(\mathbf V)=X_{f_{1,\infty}}(\mathbf U(\mathbf V)).
\]
Thus the collection 
\[
\mathcal G_{mN,f_{1,\infty}}^{*}
   :=\left\{\mathbf U(\mathbf V):\mathbf V\in\mathcal G_{N,f_{1,\infty}^m}\right\}
\]
is a cover of \(Z\) with respect to \(f_{1,\infty}\),
and
\[
m(\mathbf U(\mathbf V))=m\,m(\mathbf V)\in[mN,mN/\theta+m).
\]
Therefore, for any \(\alpha\ge0\),
\begin{equation}\label{eq:s4-en}
 \sum_{\mathbf V\in\mathcal G_{N,f_{1,\infty}^m}} e^{-\alpha m\,m(\mathbf V)}
   = \sum_{\mathbf U\in\mathcal G_{mN,f_{1,\infty}}^{*}} e^{-\alpha m(\mathbf U)}.
\end{equation}

\emph{Step~2.}
If \(m(\mathbf U(\mathbf V))<mN/\theta+1\), let 
\(\mathbf W(\mathbf V):=\mathbf U(\mathbf V)\);
otherwise, truncate further by setting
\[
q=\Big\lfloor\frac{mN}{\theta}\Big\rfloor,
\qquad
\mathbf W(\mathbf V):=\mathbf U(\mathbf V)|_{[0,q-1]}.
\]
Then \(m(\mathbf W(\mathbf V))\in[mN,mN/\theta+1)\) and
\(X_{f_{1,\infty}}(\mathbf U(\mathbf V))\subseteq X_{f_{1,\infty}}(\mathbf W(\mathbf V))\).
Thus the collection 
\[
\mathcal G_{mN,f_{1,\infty}}
   :=\left\{\mathbf W(\mathbf V):\mathbf U(\mathbf V) \in \mathcal G_{mN,f_{1,\infty}}^{*}\right\}
\]
is a cover of \(Z\) with respect to \(f_{1,\infty}\).
When the additional truncation is applied, we have 
\[
0
\le m(\mathbf U(\mathbf V))-m(\mathbf W(\mathbf V))
< \frac{mN}{\theta}+m-\Bigl\lfloor\frac{mN}{\theta}\Bigr\rfloor
< m+1.
\]
Therefore,
\begin{equation}\label{eq:s5-en}
\sum_{\mathbf U\in\mathcal G_{mN,f_{1,\infty}}^{*}}e^{-\alpha m(\mathbf U)}
   \ge e^{-\alpha(m+1)}
      \sum_{\mathbf W\in\mathcal G_{mN,f_{1,\infty}}}e^{-\alpha m(\mathbf W)}.
\end{equation}

\emph{Step~3.}
For each \(\mathbf W \in\mathcal G_{mN,f_{1,\infty}}\) and every tuple 
\((U_{1},\ldots,U_{r}) \in \mathcal{U}^{r}\), define
\[
\mathbf U^*(\mathbf W) := (\mathbf W, U_{1},\ldots,U_{r}),
\]
so that \(m(\mathbf U^*(\mathbf W)) \in [k,k/\theta+1)\).
Since
\[
\bigcup_{(U_1,\ldots,U_r)\in\mathcal U^r}
   X_{f_{1,\infty}}\!\bigl(\mathbf U^*(\mathbf W)\bigr)
   = X_{f_{1,\infty}}\!\bigl(\mathbf W\bigr),
\]
the collection
\[
\mathcal G_{k,f_{1,\infty}}
   := \left\{\mathbf U^*(\mathbf W) :
       \mathbf W \in\mathcal G_{mN,f_{1,\infty}}\right\}
\]
still covers \(Z\). Therefore,
\begin{equation}\label{eq:s6-en}
\mathrm{card}(\mathcal{U}^{r})
   \sum_{\mathbf W \in \mathcal G_{mN,f_{1,\infty}}}
   e^{-\alpha m(\mathbf W)}
   = e^{\alpha r}\!\!\sum_{\mathbf U^*\in\mathcal G_{k,f_{1,\infty}}}
     e^{-\alpha m(\mathbf U^*)}.
\end{equation}

Combining \eqref{eq:s4-en}–\eqref{eq:s6-en}, we obtain
\[
(\mathrm{card}(\mathcal{U}))^{r} e^{\alpha(m+1-r)}
   \sum_{\mathbf V \in \mathcal G_{N,f_{1,\infty}^m}}
   e^{-\alpha m\,m(\mathbf V)}
   \ge 
   \sum_{\mathbf U^*\in\mathcal G_{k,f_{1,\infty}}}
   e^{-\alpha m(\mathbf U^*)},
\]
which implies
\[
(\mathrm{card}(\mathcal{U}))^{m} e^{\alpha(m+1)} M_{f_{1,\infty}^m}(Z,\alpha m,\mathcal V,N,\theta) \ge M_{f_{1,\infty}}\!\left(Z,\alpha,\mathcal U,k,\theta \right).
\]

Taking $\liminf_{N\to\infty}$ (resp. $\limsup_{N\to\infty}$) with $k=mN+r$, we obtain
\[
(\mathrm{card}(\mathcal{U}))^{m} e^{\alpha(m+1)} \underline{m}_{f_{1,\infty}^m}(Z,\alpha m,\mathcal V,\theta) \ge \underline{m}_{f_{1,\infty}}\!\left(Z,\alpha,\mathcal U,\theta \right),
\]
\[
(\mathrm{card}(\mathcal{U}))^{m} e^{\alpha(m+1)} \overline{m}_{f_{1,\infty}^m}(Z,\alpha m,\mathcal V,\theta) \ge \overline{m}_{f_{1,\infty}}\!\left(Z,\alpha,\mathcal U,\theta \right).
\]

It follows that
\[
\underline{h}_{\mathrm{top}}(f_{1,\infty}^m,Z,\theta)\;\ge\;
m\,\underline{h}_{\mathrm{top}}(f_{1,\infty},Z,\theta), \qquad
\overline{h}_{\mathrm{top}}(f_{1,\infty}^m,Z,\theta)\;\ge\;
m\,\overline{h}_{\mathrm{top}}(f_{1,\infty},Z,\theta),
\]
which completes the proof.
\end{proof}

\begin{proposition}\label{pr2-theta}
Let $ (X,d) $ be a compact metric space and $  f_{1,\infty} $ be a sequence of equicontinuous maps from X to itself. Then for any $Z\subseteq X$, $\theta \in [0,1]$ and $m \in \mathbb{N}$,
\[
\underline{h}_{\mathrm{top}}(f_{1,\infty}^m,Z,\theta)=
m\,\underline{h}_{\mathrm{top}}(f_{1,\infty},Z,\theta),
\qquad
\overline{h}_{\mathrm{top}}(f_{1,\infty}^m,Z,\theta)=
m\,\overline{h}_{\mathrm{top}}(f_{1,\infty},Z,\theta).
\]
\end{proposition}
\begin{proof}
The inequalities
\[
\underline{h}_{\mathrm{top}}(f_{1,\infty}^m,Z,\theta)
\le m\,\underline{h}_{\mathrm{top}}(f_{1,\infty},Z,\theta)
\quad\text{and}\quad
\overline{h}_{\mathrm{top}}(f_{1,\infty}^m,Z,\theta)
\le m\,\overline{h}_{\mathrm{top}}(f_{1,\infty},Z,\theta)
\]
follow directly from Proposition~\ref{pr-theta}.
It remains to establish the reverse inequalities.

Fix $\varepsilon>0$ and define
\[
\delta(\varepsilon)
=\varepsilon+\sup_{i\ge1}\max_{1\le k\le m-1}
\sup_{d(x,y)\le\varepsilon} d\big(f_i^k(x),f_i^k(y)\big),
\]
so that $\delta(\varepsilon)\to0$ as $\varepsilon\to0$.  
For any $N\in\mathbb N$, let 
\[
\mathcal G_{N,f_{1,\infty}^m}
=\left\{B_{n_i,f_{1,\infty}^m}(x_i,\varepsilon)\right\}_i
\]
be a cover of $Z$, where each $B_{n_i,f_{1,\infty}^m}(x_i,\varepsilon)$ is the $(n_i,\varepsilon)$–Bowen ball 
with respect to $f_{1,\infty}^m$ and $n_i\in[N,\,N/\theta+1)$. 
By equicontinuity, for each $i$,
\[
B_{n_i,f_{1,\infty}^m}(x_i,\varepsilon)
\subseteq B_{m n_i,f_{1,\infty}}\left(x_i,\delta(\varepsilon)\right).
\]
Hence the family
\[
\mathcal H
=\Bigl\{
B_{m n_i,f_{1,\infty}}(x_i,\delta(\varepsilon)):
B_{n_i,f_{1,\infty}^m}(x_i,\varepsilon)
\in\mathcal G_{N,f_{1,\infty}^m}
\Bigr\}
\]
is a cover of $Z$ by $(mn_i,\delta(\varepsilon))$–Bowen balls with respect to $f_{1,\infty}$ with lengths in $[mN,\,mN/\theta+m)$. 
Moreover, for any $\alpha\ge0$,
\begin{equation}\label{eq:step1}
\sum_{B_{n_i,f_{1,\infty}^m}\in\mathcal G_{N,f_{1,\infty}^m}}
e^{-(\alpha m)n_i}
=\sum_{B_{m n_i,f_{1,\infty}}\in\mathcal H}
e^{-\alpha\,m n_i}.
\end{equation}

For each $i$, set 
\[
q_i=\min\left\{mn_i,\Bigl\lfloor\dfrac{mN}{\theta}\Bigr\rfloor\right\},
\]
then $q_i\in[mN,mN/\theta+1)$ and 
$B_{m n_i,f_{1,\infty}}(x_i,\delta(\varepsilon))
\subseteq B_{q_i,f_{1,\infty}}(x_i,\delta(\varepsilon))$.
Let 
\[
\mathcal H^*=\left\{B_{q_i,f_{1,\infty}}(x_i,\delta(\varepsilon))\right\}_i,
\]
which still covers $Z$, and since $m n_i-q_i< m+1$,
\begin{equation}\label{eq:step2}
\sum_{B_{m n_i,f_{1,\infty}}\in\mathcal H}
e^{-\alpha\,m n_i}
\ge e^{-\alpha(m+1)}
\sum_{B_{q_i,f_{1,\infty}}\in\mathcal H^*}
e^{-\alpha q_i}.
\end{equation}

Now fix $r\in\{0,\dots,m-1\}$ and put $K=mN+r$. 
We enlarge the length of each ball 
$B_{q_i,f_{1,\infty}}(x_i,\delta(\varepsilon))\in\mathcal H^*$
from $q_i$ to $q_i+r$ without losing coverage.

Since $X$ is compact, there exist points $z_1,\ldots,z_M\in X$ such that
\[
X=\bigcup_{j=1}^{M} B_d\big(z_j,\delta(\varepsilon)\big).
\]
For each $i$ and for any $u\in B_{q_i,f_{1,\infty}}(x_i,\delta(\varepsilon))$,
choose $j\in\{1,\ldots,M\}$ with 
$f_1^{q_i}(u)\in B_d(z_j,\delta(\varepsilon))$, and then pick
\[
y_{i+1,j}\in 
B_{q_i,f_{1,\infty}}(x_i,\delta(\varepsilon))\cap
f_1^{-q_i}B_d\big(z_j,\delta(\varepsilon)\big),
\]
whenever this intersection is nonempty.
By the triangle inequality, for $0\le k\le q_i-1$,
\[
d\big(f_1^k(u),f_1^k(y_{i+1,j})\big)
\le d\big(f_1^k(u),f_1^k(x_i)\big)
   + d\big(f_1^k(x_i),f_1^k(y_{i+1,j})\big)
<2\delta(\varepsilon),
\]
and for $k=q_i$,
\[
d\big(f_1^{q_i}(u),f_1^{q_i}(y_{i+1,j})\big)
\le d\big(f_1^{q_i}(u),z_j\big)
   + d\big(z_j,f_1^{q_i}(y_{i+1,j})\big)
<2\delta(\varepsilon).
\]
Hence \(u\in B_{q_i+1,f_{1,\infty}}\left(y_{i+1,j},2\delta(\varepsilon)\right)\),
and therefore
\[
B_{q_i,f_{1,\infty}}\left(x_i,\delta(\varepsilon)\right)
\subseteq 
\bigcup_{j=1}^{M}
B_{q_i+1,f_{1,\infty}}\!\left(y_{i+1,j},2\delta(\varepsilon)\right).
\]
By repeating this construction $r$ times, we obtain for each $i$ a finite family
\[
\mathcal F_i=\bigl\{
B_{q_i+r,f_{1,\infty}}\left(y_{i+r,j},2^r\delta(\varepsilon)\right)
:\ 1\le j\le M^r
\bigr\},
\]
such that
\[
B_{q_i,f_{1,\infty}}\left(x_i,\delta(\varepsilon)\right)
\subseteq
\bigcup_{B\in\mathcal F_i} B.
\]
Consequently,
\begin{equation}\label{eq:step3}
\sum_{B_{q_i,f_{1,\infty}}\in\mathcal H^*} e^{-\alpha q_i} 
\ge M^{-r} e^{\alpha r}
  \sum_i \sum_{B\in\mathcal F_i} e^{-\alpha(q_i+r)}.
\end{equation}

As $\mathcal F=\bigcup_i\mathcal F_i$ covers $Z$ and its elements have lengths in $[K,\,K/\theta+1)$, combining \eqref{eq:step1}- \eqref{eq:step3} and taking the infimum over all initial covers $\mathcal G_{N,f_{1,\infty}^m}$, we obtain
\[
M_{f_{1,\infty}^m}(Z,\alpha m,\varepsilon,N,\theta)
\ge
M^{-r}\, e^{-\alpha(m+1-r)}\,
M_{f_{1,\infty}}\left(Z,\alpha,2^r\delta(\varepsilon),K,\theta\right),
\]
which implies
\[
M^{m}\, e^{\alpha(m+1)}\,M_{f_{1,\infty}^m}\left(Z,\alpha m,\varepsilon,N,\theta\right)
\ge
M_{f_{1,\infty}}\left(Z,\alpha,2^m\delta(\varepsilon),K,\theta\right).
\]
Taking $\liminf_{N\to\infty}$ (resp. $\limsup_{N\to\infty}$) with 
$K=mN+r$, and then letting $\varepsilon\to0$ (so that $\delta(\varepsilon)\to0$),
we obtain
\[
\underline h_{\mathrm{top}}(f_{1,\infty}^m,Z,\theta)
\ge
m\,\underline h_{\mathrm{top}}(f_{1,\infty},Z,\theta),
\qquad
\overline h_{\mathrm{top}}(f_{1,\infty}^m,Z,\theta)
\ge
m\,\overline h_{\mathrm{top}}(f_{1,\infty},Z,\theta).
\]
This completes the proof.
\end{proof}

\begin{proposition}\label{le5.1-theta-refined}
Let \( X \) be a compact topological space and \( f_{1,\infty} \) be a sequence of continuous self-maps of \( X \).
For any $Z\subseteq X, k\in\mathbb N$ and $\theta\in [0,1]$, we have
\[
\underline{h}_{\mathrm{top}}\bigl(f_{k,\infty},Z,\theta\bigr)
=\underline{h}_{\mathrm{top}}\bigl(f_{k+1,\infty},f_k(Z),\theta\bigr),
\qquad
\overline{h}_{\mathrm{top}}\bigl(f_{k,\infty},Z,\theta\bigr)
=\overline{h}_{\mathrm{top}}\bigl(f_{k+1,\infty},f_k(Z),\theta\bigr),
\]
where \( f_{k,\infty}:=\left\{f_k,f_{k+1},\ldots\right\}\).
\end{proposition}
\begin{proof}
We only give the proof for the lower $\theta$-intermediate entropy. 
The argument for the upper $\theta$-intermediate entropy is entirely analogous. The proof is divided into two parts.

\medskip
\noindent\textbf{Part I.}
We prove that
\[
\underline{h}_{\rm top}(f_{k,\infty},Z,\theta)\;\ge\;
\underline{h}_{\rm top}(f_{k+1,\infty},f_k(Z),\theta).
\]
The case $\theta=0$ follows from \cite[Lemma~3.3]{ju2021pesin}, hence we assume $\theta\in(0,1]$.
Fix a finite open cover $\mathcal U$ of $X$ and $\alpha\ge0$.
For $N\in\mathbb N$, choose a family
\[
\mathcal G_{N,f_{k,\infty}}\subseteq
\bigcup_{N\le m<N/\theta+1}\mathcal S_m(\mathcal U)
\]
that covers $Z$ with respect to $f_{k,\infty}$.
Let $L:=(N-1)/\theta$ and split \(\mathcal G_{N,f_{k,\infty}}=\mathcal G_a\cup\mathcal G_b\), where
\[
\mathcal G_a:=\left\{\mathbf U:\ N\le m(\mathbf U)<L+2\right\},\quad
\mathcal G_b:=\left\{\mathbf U:\ L+2\le m(\mathbf U)<L+1/\theta+1\right\}.
\]
For $\mathbf U\in\mathcal G_a$ set $\mathbf U^{*}:=\mathbf U$; for $\mathbf U\in\mathcal G_b$ set
\(\mathbf U^{*}:=\mathbf U|_{[0,\lfloor L\rfloor]}\).
Define
\[
\widetilde{\mathcal G}:=\left\{\mathbf U^{*}:\ \mathbf U\in\mathcal G_{N,f_{k,\infty}}\right\}.
\] 
Then $\widetilde{\mathcal G}$ still covers $Z$ and \(m(\mathbf U^{*})\in [N,L+2)\) for all \(\mathbf U^{*} \in \widetilde{\mathcal G}\).
For each $\mathbf U^{*}\in\widetilde{\mathcal G}$, define its truncated string by
\[
\mathbf V(\mathbf U^{*}):=\mathbf U^{*}|_{[1,\,m(\mathbf U^{*})-1]}.
\]
Then \(m(\mathbf V(\mathbf U^{*}))\in[N-1,\,L+1)\) and
\[
f_k\!\big(X_{f_{k,\infty}}(\mathbf U^{*})\big)\subseteq
X_{f_{k+1,\infty}}\!\big(\mathbf V(\mathbf U^{*})\big),
\]
hence $\left\{\mathbf V(\mathbf U^{*}):\mathbf U^{*} \in \widetilde{\mathcal G}\right\}$ covers $f_k(Z)$ with admissible lengths.

If $\mathbf U\in\mathcal G_b$, then
\(m(\mathbf U)-m(\mathbf U^{*})< 1/\theta+1\) and thus
\[
e^{-\alpha m(\mathbf U)}\;>\;e^{-\alpha(\frac1\theta+1)}\,e^{-\alpha m(\mathbf U^{*})}
\;=\;e^{-\alpha(\frac1\theta+2)}\,e^{-\alpha m(\mathbf V(\mathbf U^{*}))}.
\]
It follows that

\begin{align*}
e^{\alpha(\frac1\theta+2)}
\sum_{\mathbf U\in\mathcal G_{N,f_{k,\infty}}}
   e^{-\alpha m(\mathbf U)}
&\;\ge\;
\sum_{\mathbf U^{*}\in\widetilde{\mathcal G}}
   e^{-\alpha m(\mathbf V(\mathbf U^{*}))}
\\[0.4em]
&\;\ge\;
M_{f_{k+1,\infty}}\bigl(f_k(Z),\alpha,\mathcal U,N-1,\theta\bigr).
\end{align*}
Taking the infimum over $\mathcal G_{N,f_{k,\infty}}$ and then $\liminf_{N\to\infty}$ gives
\[
e^{\alpha(\frac1\theta+2)}\,
\underline m_{f_{k,\infty}}(Z,\alpha,\mathcal U,\theta)
\;\ge\;
\underline m_{f_{k+1,\infty}}(f_k(Z),\alpha,\mathcal U,\theta),
\]
hence
\[
\underline h_{\mathrm{top}}(f_{k,\infty},Z,\mathcal U,\theta)
\;\ge\;
\underline h_{\mathrm{top}}(f_{k+1,\infty},f_k(Z),\mathcal U,\theta).
\]
Finally, taking the supremum over all finite open covers yields the desired inequality.

\medskip
\noindent\textbf{Part~II.}
We prove the reverse inequality
\[
\underline{h}_{\mathrm{top}}(f_{k,\infty},Z,\theta)
\;\le\;
\underline{h}_{\mathrm{top}}(f_{k+1,\infty},f_k(Z),\theta).
\]
Let $\mathcal U=\{U_1,\ldots,U_n\}$ be a finite open cover of $X$ and fix $\alpha\ge0$.
Choose
\[
\mathcal G_{N,f_{k+1,\infty}}
   \subseteq
   \bigcup_{N\le m<N/\theta+1}\mathcal S_m(\mathcal U)
\]
covering $f_k(Z)$ with respect to $f_{k+1,\infty}$.
Set
\[
\mathcal G_{N+1,f_{k,\infty}}
   :=\left\{\mathbf V(U_j,\mathbf U):
      U_j\in\mathcal U,\;
      \mathbf U\in\mathcal G_{N,f_{k+1,\infty}}\right\}.
\]
Then for each $\mathbf V(U_j,\mathbf U)\in\mathcal G_{N+1,f_{k,\infty}}$,
\[
X_{f_{k,\infty}}(\mathbf V(U_j,\mathbf U))
   =U_j\cap f_k^{-1}\bigl(X_{f_{k+1,\infty}}(\mathbf U)\bigr).
\]
It follows that
\begin{align*}
Z \subseteq f_k^{-1}(f_k(Z))
     &\subseteq 
     \bigcup_{U_j\in\mathcal U} \bigcup_{\mathbf U\in\mathcal G_{N,f_{k+1,\infty}}}
     U_j \cap f_k^{-1}\!\bigl(X_{f_{k+1,\infty}}(\mathbf U)\bigr) \\
   &=\bigcup_{\mathbf V(U_j,\mathbf U)\in\mathcal G_{N+1,f_{k,\infty}}}
     X_{f_{k,\infty}}(\mathbf V(U_j,\mathbf U)).   
\end{align*}
Thus $\mathcal G_{N+1,f_{k,\infty}}$ covers $Z$ with respect to $f_{k,\infty}$ and $m(\mathbf V)\in[N+1,(N+1)/\theta+1)$ for all \(\mathbf V \in \mathcal G_{N+1,f_{k,\infty}}\). Then
\begin{align*}
n\, e^{-\alpha}
   \sum_{\mathbf U\in\mathcal G_{N,f_{k+1,\infty}}}
      e^{-\alpha m(\mathbf U)}
&=\sum_{U_j\in\mathcal U}
   \sum_{\mathbf U\in\mathcal G_{N,f_{k+1,\infty}}}
      e^{-\alpha(m(\mathbf U)+1)}\\[0.3em]
&=\sum_{\mathbf V(U_j,\mathbf U)\in\mathcal G_{N+1,f_{k,\infty}}}
      e^{-\alpha m(\mathbf V(U_j,\mathbf U))}\\[0.3em]
&\ge M_{f_{k,\infty}}(Z,\alpha,\mathcal U,N+1,\theta).
\end{align*}
Taking the infimum over $\mathcal G_{N,f_{k+1,\infty}}$ and then $\liminf_{N\to\infty}$ gives
\[
\underline m_{f_{k,\infty}}(Z,\alpha,\mathcal U,\theta)\le 
n\,e^{-\alpha}\,
\underline m_{f_{k+1,\infty}}(f_{k}(Z),\alpha,\mathcal U,\theta),
\]
which implies
\[
\underline h_{\rm top}(f_{k,\infty},Z,\mathcal{U},\theta)\le \underline h_{\rm top}(f_{k+1,\infty},f_{k}(Z),\mathcal{U},\theta).
\]
Taking the supremum over all finite open covers yields 
\[
\underline h_{\mathrm{top}}(f_{k,\infty},Z,\theta)\le \underline h_{\rm top}(f_{k+1,\infty},f_{k}(Z),\theta).
\]

Combining Parts~I and~II completes the proof.
\end{proof}

\begin{corollary}\label{momo-theta}
Let $X$ be a compact topological space and $Z \subseteq X$.  
We say that $Z$ is:
\begin{itemize}
\item \emph{$f_{1,\infty}$-forward invariant} if $f_k(Z)\subseteq Z$ for all $k\in\mathbb N$;
\item \emph{$f_{1,\infty}$-backward invariant} if $Z\subseteq f_k(Z)$ for all $k\in\mathbb N$;
\item \emph{$f_{1,\infty}$-invariant} if $f_k(Z)=Z$ for all $k\in\mathbb N$.
\end{itemize}

Then for any \(1 \le i < j < \infty\) and $\theta \in [0,1]$, the following hold:
\begin{enumerate}
\item[(1)] If $Z$ is forward invariant, then
\[
\underline{h}_{\mathrm{top}}(f_{i,\infty},Z,\theta)
\le
\underline{h}_{\mathrm{top}}(f_{j,\infty},Z,\theta),
\qquad
\overline{h}_{\mathrm{top}}(f_{i,\infty},Z,\theta)
\le
\overline{h}_{\mathrm{top}}(f_{j,\infty},Z,\theta).
\]

\item[(2)] If $Z$ is backward invariant, then
\[
\underline{h}_{\mathrm{top}}(f_{i,\infty},Z,\theta)
\ge
\underline{h}_{\mathrm{top}}(f_{j,\infty},Z,\theta),
\qquad
\overline{h}_{\mathrm{top}}(f_{i,\infty},Z,\theta)
\ge
\overline{h}_{\mathrm{top}}(f_{j,\infty},Z,\theta).
\]

\item[(3)] If $Z$ is invariant, then
\[
\underline{h}_{\mathrm{top}}(f_{i,\infty},Z,\theta)
=\underline{h}_{\mathrm{top}}(f_{j,\infty},Z,\theta),
\qquad
\overline{h}_{\mathrm{top}}(f_{i,\infty},Z,\theta)
=\overline{h}_{\mathrm{top}}(f_{j,\infty},Z,\theta).
\]
\end{enumerate}
\end{corollary}
\begin{proof}
By Proposition~\ref{le5.1-theta-refined}, for every $k \in \mathbb{N}$,
\[
\underline h_{\mathrm{top}}(f_{k,\infty},Z,\theta)
=\underline h_{\mathrm{top}}(f_{k+1,\infty},f_k(Z),\theta),
\quad
\overline h_{\mathrm{top}}(f_{k,\infty},Z,\theta)
=\overline h_{\mathrm{top}}(f_{k+1,\infty},f_k(Z),\theta).
\]
Iterating from $k=i$ to $j-1$ yields
\[
\underline h_{\mathrm{top}}(f_{i,\infty},Z,\theta)
=\underline h_{\mathrm{top}}\bigl(f_{j,\infty},\,f_i^{j-i}(Z),\theta\bigr),
\]
and the analogous identity for $\overline h_{\mathrm{top}}$.
If $Z$ is forward invariant, then $f_i^{j-i}(Z) \subseteq Z$, hence
\(\underline h_{\mathrm{top}}(f_{i,\infty},Z,\theta)\le \underline h_{\mathrm{top}}(f_{j,\infty},Z,\theta)\)
(and similarly for $\overline h_{\mathrm{top}}$).
If $Z$ is backward invariant, the inclusion reverses and so do the inequalities.
If $Z$ is invariant, equalities hold.
\end{proof}

\begin{proposition}\label{commute-theta}
Let $f_1,f_2$ be continuous self-maps on a compact topological space $X$ and let $Z\subseteq X$.
If $Z$ is $\{f_1,f_2\}$-forward invariant or $\{f_1,f_2\}$-backward invariant, then for every $\theta\in[0,1]$ we have
\[
\underline{h}_{\mathrm{top}}(f_1\!\circ\! f_2,Z,\theta)
=\underline{h}_{\mathrm{top}}(f_2\!\circ\! f_1,Z,\theta),
\qquad
\overline{h}_{\mathrm{top}}(f_1\!\circ\! f_2,Z,\theta)
=\overline{h}_{\mathrm{top}}(f_2\!\circ\! f_1,Z,\theta).
\]
\end{proposition}
\begin{proof}
Let
\[
f_{1,\infty}=\left\{f_1,f_2,f_1,f_2,\ldots\right\},
\qquad
g_{1,\infty}=\left\{f_2,f_1,f_2,f_1,\ldots\right\}.
\]
By Corollary~\ref{momo-theta}, if $Z$ is $\{f_1,f_2\}$-forward invariant, then
\[
\underline{h}_{\mathrm{top}}(f_{1,\infty},Z,\theta)
\le
\underline{h}_{\mathrm{top}}(g_{1,\infty},Z,\theta)
\le
\underline{h}_{\mathrm{top}}(f_{1,\infty},Z,\theta),
\]
\[
\overline{h}_{\mathrm{top}}(f_{1,\infty},Z,\theta)
\le
\overline{h}_{\mathrm{top}}(g_{1,\infty},Z,\theta)
\le
\overline{h}_{\mathrm{top}}(f_{1,\infty},Z,\theta),
\]
and the inequalities reverse for backward invariance.  
Hence
\[
\underline{h}_{\mathrm{top}}(f_{1,\infty},Z,\theta)
=\underline{h}_{\mathrm{top}}(g_{1,\infty},Z,\theta),
\qquad
\overline{h}_{\mathrm{top}}(f_{1,\infty},Z,\theta)
=\overline{h}_{\mathrm{top}}(g_{1,\infty},Z,\theta).
\]
Applying Proposition~\ref{pr1-theta} with $m=2$,
we obtain
\[
\underline{h}_{\mathrm{top}}\left(f_1\!\circ\! f_2,Z,\theta \right)
=\underline{h}_{\mathrm{top}}(f_2\!\circ\! f_1,Z,\theta),
\qquad
\overline{h}_{\mathrm{top}}(f_1\!\circ\! f_2,Z,\theta)
=\overline{h}_{\mathrm{top}}(f_2\!\circ\! f_1,Z,\theta).
\]
\end{proof}

\begin{corollary}\label{cycle-theta}
Let $X$ be a compact topological space, and let $f_{1},\dots,f_{n}$ be continuous self-maps on $X$. 
Suppose $Z\subseteq X$ is $\{f_{1},\dots,f_{n}\}$-forward invariant or $\{f_{1},\dots,f_{n}\}$-backward invariant. 
Then for any $1<i\le n$ and $\theta\in[0,1]$,
\[
\underline{h}_{\mathrm{top}}\left(f_{n}\circ\cdots\circ f_{2}\circ f_{1},Z,\theta\right)
=\underline{h}_{\mathrm{top}}\left(f_{i-1}\circ\cdots\circ f_{2}\circ f_{1}\circ f_{n}\circ\cdots\circ f_{i},Z,\theta\right),
\]
\[
\overline{h}_{\mathrm{top}}\left(f_{n}\circ\cdots\circ f_{2}\circ f_{1},Z,\theta\right)
=\overline{h}_{\mathrm{top}}\left(f_{i-1}\circ\cdots\circ f_{2}\circ f_{1}\circ f_{n}\circ\cdots\circ f_{i},Z,\theta\right).
\]
\end{corollary}
\begin{proof}
Let 
\[
g=f_{n}\circ\cdots\circ f_{i},\qquad h=f_{i-1}\circ\cdots\circ f_{1}.
\]
Then, by Proposition~\ref{commute-theta},
\[
\underline{h}_{\mathrm{top}}(g\circ h,Z,\theta)
=\underline{h}_{\mathrm{top}}(h\circ g,Z,\theta),
\qquad
\overline{h}_{\mathrm{top}}(g\circ h,Z,\theta)
=\overline{h}_{\mathrm{top}}(h\circ g,Z,\theta),
\]
which yields the claim.
\end{proof}

\begin{proposition}\label{prop:product-inequality}
Let \((X, d)\) and \((Y, \rho)\) be compact metric spaces, and let
$f_{1,\infty}$ and $g_{1,\infty}$ be sequences of continuous self-maps on $X$ and $Y$, respectively. Define a metric \(\mathbf d\) on \(X \times Y\) by \(\mathbf d((x_1, y_1), (x_2, y_2)) = \max \{ d(x_1, x_2), \rho(y_1, y_2) \}\) and a sequence of transformations on \(X \times Y\) by \((f_i \times g_i)(x, y) = (f_i x, g_i y)\). Then for all nonempty $Z\subseteq X$, $W\subseteq Y$ and every $\theta\in[0,1]$, we have
\[
\underline{h}_{\mathrm{top}}(f_{1,\infty}\times g_{1,\infty},Z \times W,\theta)
 \le \underline{h}_{\mathrm{top}}(f_{1,\infty},Z,\theta)+h(g_{1,\infty};W),
\]
\[
\overline{h}_{\mathrm{top}}(f_{1,\infty}\times g_{1,\infty},Z \times W,\theta)
 \le \overline{h}_{\mathrm{top}}(f_{1,\infty},Z,\theta)+h(g_{1,\infty};W).
\] 
\end{proposition}
\begin{proof}
For clarity, we denote by $B^X_n(x,\xi)$, $B^Y_n(y,\xi)$ and 
$B^{X\times Y}_n((x,y),\xi)$ the $(n,\xi)$–Bowen balls with respect to
$f_{1,\infty}$ on $X$, $g_{1,\infty}$ on $Y$ and
$f_{1,\infty}\times g_{1,\infty}$ on $X\times Y$, respectively.
Let $t>h(g_{1,\infty};W)$ and $s>\overline h_{\mathrm{top}}(f_{1,\infty},Z,\theta)$.
For any $\xi>0$, there exists $N_0\in\mathbb N$ such that
\[
r_n(g_{1,\infty},W,\xi)\le e^{tn}\qquad\text{for all }n\ge N_0.
\]
Moreover, for any $\zeta>0$ there exists $N_1\in\mathbb N$ such that for all $N\ge N_1$
there is a finite (or countable) family $\{B^X_{n_i}(x_i,\xi)\}_i$ satisfying
\[
N\le n_i<N/\theta+1,\qquad 
Z\subseteq \bigcup_i B^X_{n_i}(x_i,\xi),\qquad 
\sum_i e^{-s n_i}<\zeta.
\]

Fix $N\ge \max\{N_0,N_1\}$. For each $i$, choose an $(n_i,\xi)$-spanning set
$E_{n_i}\subseteq W$ with $\mathrm{card}(E_{n_i})=r_{n_i}(g_{1,\infty},W,\xi)$.
Since
\[
B^{X\times Y}_{n_i}\bigl((x_i,y),\xi\bigr)
= B^X_{n_i}(x_i,\xi)\times B^Y_{n_i}(y,\xi),
\]
then the family $\left\{B^{X\times Y}_{n_i}((x_i,y),\xi):y\in E_{n_i}\right\}_i$ covers $Z\times W$.
Thus,
\[
\begin{aligned}
M(Z\times W,s+t,N,\xi,\theta)
&\le \sum_i\sum_{y\in E_{n_i}} e^{-(s+t)n_i}
= \sum_i r_{n_i}(g_{1,\infty},W,\xi)\,e^{-(s+t)n_i}\\
&\le \sum_i e^{-s n_i}
< \zeta.
\end{aligned}
\]
As $\zeta>0$ is arbitrary, $\overline m(Z\times W,s+t,\xi,\theta)=0$,
and therefore
\[
\overline h_{\mathrm{top}}(f_{1,\infty}\times g_{1,\infty},Z\times W,\theta)
\le s+t.
\]
Letting $s\downarrow \overline h_{\mathrm{top}}(f_{1,\infty},Z,\theta)$
and $t\downarrow h(g_{1,\infty};W)$ yields the claimed inequality.
The argument with $\liminf$ instead of $\limsup$ gives the bound for 
$\underline h_{\mathrm{top}}$.
\end{proof}

\begin{remark}
When $\theta=0$, we have
\[
h_{\mathrm{top}}^{B}(f_{1,\infty}\times g_{1,\infty}, Z \times W)
\le 
h_{\mathrm{top}}^{B}(f_{1,\infty}, Z)
+
h_{\mathrm{top}}^{P}(g_{1,\infty}, W),
\]
as shown in \cite[Theorem~6.1]{chen2025nonautonomous}, and
\[
h_{\mathrm{top}}^{B}(g_{1,\infty}, W)
\le 
h_{\mathrm{top}}^{P}(g_{1,\infty}, W)
\le 
h(g_{1,\infty};W),
\]
see \cite[Proposition~4.7]{chen2025nonautonomous}.
Here $h_{\mathrm{top}}^{P}(g_{1,\infty}, W)$ denotes the packing topological entropy on the set $W$, which may be strictly smaller than the upper capacity entropy. Therefore, the inequalities in Proposition~\ref{prop:product-inequality} can be strict.
\end{remark}

\section{Topological conjugacy}\label{sec4}
In this section, we investigate the relations between the $\theta$-intermediate topological entropies of two topologically equisemiconjugate nonautonomous dynamical systems.

Let $(X, f_{1,\infty})$ and $(Y, g_{1,\infty})$ be two NDSs,
where $f_{1,\infty}=\{f_i:X\to X\}_{i=1}^{\infty}$ and $g_{1,\infty}=\{g_i:Y\to Y\}_{i=1}^{\infty}$ are sequences of continuous maps.

A sequence $\pi_{1,\infty}=\{\pi_i\}_{i=1}^{\infty}$ of surjective continuous maps 
$\pi_i:X\to Y$ is called a sequence of factor maps
or a semiconjugacy from $(X,f_{1,\infty})$ to $(Y,g_{1,\infty})$ if the following commutation relation holds:
\[
\pi_{i+1}\circ f_i = g_i\circ \pi_i \qquad\text{for all } i\ge1.
\]

If such a sequence exists, we say that $(Y,g_{1,\infty})$ is a \emph{factor} of $(X,f_{1,\infty})$ or $(X,f_{1,\infty})$ is an \emph{extension} of $(Y,g_{1,\infty})$.

\begin{definition}
Let $\pi_{1,\infty}=\{\pi_i\}_{i=1}^{\infty}$ be a semiconjugacy from $(X,f_{1,\infty})$ to $(Y,g_{1,\infty})$. If $\pi_{1,\infty}$ is equicontinuous, we call it a 
(topological) equisemiconjugacy, and say that $(X,f_{1,\infty})$ is
(topologically) equisemiconjugate to $(Y,g_{1,\infty})$.
Furthermore, if each $\pi_i$ is a homeomorphism and the inverse sequence 
$\pi_{1,\infty}^{-1}=\{\pi_i^{-1}\}_{i=1}^{\infty}$ is also equicontinuous, then
$\pi_{1,\infty}$ is called a (topological) equiconjugacy, and the two systems are (topologically) equiconjugate.    
\end{definition}

\begin{theorem}\label{con}
Let $(X,d)$ and $(Y,\rho)$ be compact metric spaces, and let 
$f_{1,\infty}$ and $g_{1,\infty}$ be sequences of continuous self-maps on $X$ and $Y$, respectively.  
Suppose $\pi_{1,\infty}=\{\pi_i\}_{i=1}^\infty$ is a topological equisemiconjugacy from $(X,f_{1,\infty})$ to $(Y,g_{1,\infty})$.  
Then for every $Z\subseteq X$, $k\in\mathbb N$, and $\theta\in[0,1]$,
\[
\underline{h}_{\mathrm{top}}(f_{k,\infty},Z,\theta)\;\ge\;\underline{h}_{\mathrm{top}}(g_{k,\infty},\pi_k(Z),\theta),
\qquad
\overline{h}_{\mathrm{top}}(f_{k,\infty},Z,\theta)\;\ge\;\overline{h}_{\mathrm{top}}(g_{k,\infty},\pi_k(Z),\theta).
\]
If $\pi_{1,\infty}$ is a topological equiconjugacy, then equality holds:
\[
\underline{h}_{\mathrm{top}}(f_{k,\infty},Z,\theta)\;=\;\underline{h}_{\mathrm{top}}(g_{k,\infty},\pi_k(Z),\theta),
\qquad
\overline{h}_{\mathrm{top}}(f_{k,\infty},Z,\theta)\;=\;\overline{h}_{\mathrm{top}}(g_{k,\infty},\pi_k(Z),\theta).
\]
\end{theorem}
\begin{proof}
Since $\{\pi_i\}_{i\ge1}$ is equicontinuous, for any $ \varepsilon>0$ there exists $0<\delta<\varepsilon$ such that 
\[
d(x,y)<\delta \quad \Longrightarrow \quad \rho(\pi_i(x),\pi_i(y))<\varepsilon,\quad \forall i\ge1.
\]
Let 
\[
\mathcal G_{N,f_{k,\infty}}
=\bigl\{\,B_{n_i,f_{k,\infty}}(x_i,\delta)\,\bigr\}_i
\]
be a cover of \(Z\), where each \(B_{n_i,f_{k,\infty}}(x_i,\delta)\) is the \((n_i,\delta)\)–Bowen ball 
with respect to \(f_{k,\infty}\) and \(n_i\in[N,\,N/\theta+1)\).
Then for every 
\(B_{n_i,f_{k,\infty}}(x_i,\delta)\in\mathcal G_{N,f_{k,\infty}}\), we have
\[
B_{n_i,f_{k,\infty}}(x_i,\delta)
=\bigcap_{p=0}^{n_i-1} f_k^{-p}\!\left(B_d(f_k^{p}(x_i),\delta)\right).
\]
For any \(y\in f_k^{-p}\!\left(B_d(f_k^{p}(x_i),\delta)\right)\) with \(0\le p\le n_i-1\),
we have \(d(f_k^{p}(x_i),f_k^{p}(y))<\delta\), and thus
\[
\rho\!\left(\pi_{k+p}f_k^{p}(x_i),\,\pi_{k+p}f_k^{p}(y)\right)<\varepsilon.
\]
Therefore
\[
y\in(\pi_{k+p}\circ f_k^{p})^{-1}
   \!\left(B_{\rho}(\pi_{k+p}f_k^{p}(x_i),\varepsilon)\right)
   =(g_k^{p}\pi_k)^{-1}\!\left(B_{\rho}(g_k^{p}\pi_k(x_i),\varepsilon)\right).
\]
Consequently,
\[
B_{n_i,f_{k,\infty}}(x_i,\delta)
\subseteq
\pi_k^{-1}\!\left(B_{n_i,g_{k,\infty}}(\pi_k(x_i),\varepsilon)\right).
\]

Hence the family
\[
\mathcal G_{N,g_{k,\infty}}
=\left\{
B_{n_i,g_{k,\infty}}(\pi_k(x_i),\varepsilon):
B_{n_i,f_{k,\infty}}(x_i,\delta)\in\mathcal G_{N,f_{k,\infty}}
\right\}
\]
covers \(\pi_k(Z)\). Moreover,
\begin{align*}
M_{f_{k,\infty}}(Z,\alpha,\delta,N,\theta)
&=\inf_{\mathcal G_{N,f_{k,\infty}}}
   \sum_{B_{n_i,f_{k,\infty}}(x_i,\delta)\in\mathcal G_{N,f_{k,\infty}}}
   e^{-\alpha n_i}\\
&\ge
\inf_{\mathcal G_{N,g_{k,\infty}}}
   \sum_{B_{n_i,g_{k,\infty}}(\pi_k(x_i),\varepsilon)\in\mathcal G_{N,g_{k,\infty}}}
   e^{-\alpha n_i}\\
&\ge
M_{g_{k,\infty}}(\pi_k(Z),\alpha,\varepsilon,N,\theta).
\end{align*}

Taking \(\liminf_{N\to\infty}\) yields
\[
\underline{m}_{f_{k,\infty}}(Z,\alpha,\delta,\theta)\;\ge\;\underline{m}_{g_{k,\infty}}(\pi_k(Z),\alpha,\varepsilon,\theta),
\]
which implies
\[
\underline{h}_{\mathrm{top}}(f_{k,\infty},Z,\delta,\theta)\;\ge\;\underline{h}_{\mathrm{top}}(g_{k,\infty},\pi_k(Z),\varepsilon,\theta).
\]

Finally, letting \(\varepsilon\to0\) gives
\[
\underline h_{\mathrm{top}}(f_{k,\infty},Z,\theta)
\;\ge\;
\underline h_{\mathrm{top}}(g_{k,\infty},\pi_k(Z),\theta).
\]
In addition, if each $\pi_i$ is a homeomorphism and $\{\pi_i^{-1}\}$ is equicontinuous, the same argument applied to $\pi_{1,\infty}^{-1}$ shows the reverse inequality.  
Hence in this case we obtain equality:
\[
\underline{h}_{\mathrm{top}}(f_{k,\infty},Z,\theta)=\underline{h}_{\mathrm{top}}(g_{k,\infty},\pi_k(Z),\theta).
\]
The proof for the upper \(\theta\)-intermediate entropy is entirely analogous.
\end{proof}

\begin{remark}
The commutative diagram
\[
\begin{tikzcd}[column sep=2.5em, row sep=2.5em]
X \arrow[r, "f_1"] \arrow[d, "f_1"] 
  & X \arrow[r, "f_2"] \arrow[d, "f_2"] 
  & \cdots \arrow[r, "f_{j-1}"] 
  & X \arrow[r, "f_j"] \arrow[d, "f_j"] 
  & X \arrow[r, "f_{j+1}"] \arrow[d, "f_{j+1}"] 
  & \cdots \\
X \arrow[r, "f_2"] 
  & X \arrow[r, "f_3"] 
  & \cdots \arrow[r, "f_j"]
  & X \arrow[r, "f_{j+1}"] 
  & X \arrow[r, "f_{j+2}"] 
  & \cdots
\end{tikzcd}
\]
shows that by setting \(\pi_i=f_i\) and \(g_i=f_{i+1}\) for all \(i\ge1\), one obtains a topological equisemiconjugacy from
\((X,f_{1,\infty})\) to \((X,f_{2,\infty})\) if $f_{1,\infty}$ is a sequence of equicontinuous surjective maps.
Fix $k\in\mathbb N$. By Theorem~\ref{con} we obtain
\[
\underline{h}_{\mathrm{top}}(f_{k,\infty},Z,\theta)
\;\ge\;
\underline{h}_{\mathrm{top}}(f_{k+1,\infty},f_k(Z),\theta),
\quad
\overline{h}_{\mathrm{top}}(f_{k,\infty},Z,\theta)
\;\ge\;
\overline{h}_{\mathrm{top}}(f_{k+1,\infty},f_k(Z),\theta).
\]
In addition, if $(X,f_{1,\infty})$ and $(X,f_{2,\infty})$ are topologically equiconjugate, then the two inequalities become equalities. This observation provides an alternative proof of Proposition~\ref{le5.1-theta-refined}, but under the stronger assumption that 
$f_{1,\infty}$ is a topological equiconjugacy between $(X,f_{1,\infty})$ and $(X,f_{2,\infty})$.
\end{remark}

\begin{corollary}
If $g:X \to X$ is a homeomorphism commuting with $f_{1,\infty}$ 
(i.e.\ $f_i\circ g=g\circ f_i$ for all $i\ge1$), then for any $Z\subseteq X$ 
and $\theta\in[0,1]$,
\[
\underline{h}_{\mathrm{top}}(f_{1,\infty},Z,\theta)
=\underline{h}_{\mathrm{top}}(f_{1,\infty},g(Z),\theta),
\qquad
\overline{h}_{\mathrm{top}}(f_{1,\infty},Z,\theta)
=\overline{h}_{\mathrm{top}}(f_{1,\infty},g(Z),\theta).
\]
\end{corollary}
\begin{proof}
The commutativity condition gives a commutative diagram for each $j$:
\[
\begin{tikzcd}[column sep=2.5em, row sep=2.5em]
X \arrow[r, "f_1"] \arrow[d, "g"] 
  & X \arrow[r, "f_2"] \arrow[d, "g"] 
  & \cdots \arrow[r, "f_{j-1}"] 
  & X \arrow[r, "f_j"] \arrow[d, "g"] 
  & X \arrow[r, "f_{j+1}"] \arrow[d, "g"] 
  & \dotsb \\
X \arrow[r, "f_1"] 
  & X \arrow[r, "f_2"] 
  & \cdots \arrow[r, "f_{j-1}"] 
  & X \arrow[r, "f_j"] 
  & X \arrow[r, "f_{j+1}"] 
  & \dotsb
\end{tikzcd}
\]
Thus the sequence $\pi_i:=g$ gives an equiconjugacy from $(X,f_{1,\infty})$ 
to itself. By Theorem~\ref{con}, equiconjugacy preserves 
both lower and upper $\theta$-intermediate topological entropies. Hence
\[
\underline{h}_{\mathrm{top}}(f_{1,\infty},Z,\theta)
=\underline{h}_{\mathrm{top}}(f_{1,\infty},g(Z),\theta),
\qquad
\overline{h}_{\mathrm{top}}(f_{1,\infty},Z,\theta)
=\overline{h}_{\mathrm{top}}(f_{1,\infty},g(Z),\theta).
\]
\end{proof}

\begin{corollary}
Let $(X,d)$ and $(Y,\rho)$ be compact metric spaces, and let
$f_{1,\infty}$ and $g_{1,\infty}$ be sequences of continuous self-maps on $X$ and $Y$, respectively.
For all nonempty subsets $Z\subseteq X$, $W\subseteq Y$, and every $\theta\in[0,1]$, we have
\[
\max\Big\{
\overline{h}_{\mathrm{top}}(f_{1,\infty},\, Z,\, \theta),\;
\overline{h}_{\mathrm{top}}(g_{1,\infty},\, W,\, \theta)
\Big\}\le \overline{h}_{\mathrm{top}}(f_{1,\infty}\times g_{1,\infty},\, Z\times W,\, \theta),
\]
\[
\max\Big\{
\underline{h}_{\mathrm{top}}(f_{1,\infty},\, Z,\, \theta),\;
\underline{h}_{\mathrm{top}}(g_{1,\infty},\, W,\, \theta)
\Big\}\le \underline{h}_{\mathrm{top}}(f_{1,\infty}\times g_{1,\infty},\, Z\times W,\, \theta).
\]
\end{corollary}
\begin{proof}
We prove the inequality for the upper intermediate entropy $\overline{h}_{\mathrm{top}}$; 
the argument for the lower entropy $\underline{h}_{\mathrm{top}}$ is analogous. 
Let $\pi_X:X\times Y\to X$ be the canonical projection and set 
$\pi_{1,\infty}=\{\pi_i\}_{i\ge1}$ with $\pi_i=\pi_X$ for all $i$. 
Each $\pi_i$ is continuous and surjective, and satisfies 
$\pi_{i+1}\circ(f_i\times g_i)=f_i\circ\pi_i$, 
so that $\pi_{1,\infty}$ is a topological equisemiconjugacy from 
$(X\times Y,f_{1,\infty}\times g_{1,\infty})$ to $(X,f_{1,\infty})$. 
By Theorem~\ref{con}, 
\[
\overline{h}_{\mathrm{top}}(f_{1,\infty}\times g_{1,\infty},Z\times W,\theta)
\ge
\overline{h}_{\mathrm{top}}(f_{1,\infty},\pi_X(Z\times W),\theta)
=
\overline{h}_{\mathrm{top}}(f_{1,\infty},Z,\theta).
\]
Applying the same argument to the projection $\pi_Y:X\times Y\to Y$ gives 
\[
\overline{h}_{\mathrm{top}}(f_{1,\infty}\times g_{1,\infty},Z\times W,\theta)
\ge
\overline{h}_{\mathrm{top}}(g_{1,\infty},W,\theta).
\]
Hence,
\[
\overline{h}_{\mathrm{top}}(f_{1,\infty}\times g_{1,\infty},Z\times W,\theta)
\ge
\max\Big\{
\overline{h}_{\mathrm{top}}(f_{1,\infty},Z,\theta),\,
\overline{h}_{\mathrm{top}}(g_{1,\infty},W,\theta)
\Big\}.
\]
\end{proof}

To establish the following theorem, we first recall the notion of \emph{topological sup-entropy} proposed by Kolyada and Snoha \cite{kolyada1996topological}.
Let $(X,d)$ be a compact metric space, $Z\subseteq X$ a nonempty subset, and
$f_{1,\infty}=\{f_i\}_{i\ge1}$ a sequence of equicontinuous self-maps of $X$.
For each $n\ge1$, define
\[
d_n^*(x,y)=\sup_i\,\max_{0\le j<n} d(f_i^j(x),f_i^j(y)), \qquad x,y\in X.
\]
Since $f_{1,\infty}$ is equicontinuous, $d_n^*$ is equivalent to $d$, and thus $(X,d_n^*)$ is also compact.

A subset $E^*\subseteq X$ is said to be $(n,\varepsilon)^*$-separated if $d_n^*(x,y)>\varepsilon$ for any distinct $x,y\in E^*$.
A set $F^*\subseteq X$ $(n,\varepsilon)^*$-spans $Z$ if for every $x\in Z$ there exists $y\in F^*$ such that $d_n^*(x,y)\le\varepsilon$.
Let \(s_n^*(f_{1,\infty};Z;\varepsilon)\) denote the maximal cardinality of an $(n,\varepsilon)^*$-separated set in $Z$, and
\(r_n^*(f_{1,\infty};Z;\varepsilon)\) the minimal cardinality of an $(n,\varepsilon)^*$-spanning set in $Z$.
The \emph{topological sup-entropy} of $f_{1,\infty}$ on $Z$ is then defined by
\[
H(f_{1,\infty};Z)
=\lim_{\varepsilon\to0}\limsup_{n\to\infty}\frac{1}{n}\log r_n^*(f_{1,\infty};Z;\varepsilon)
=\lim_{\varepsilon\to0}\limsup_{n\to\infty}\frac{1}{n}\log s_n^*(f_{1,\infty};Z;\varepsilon).
\]

\begin{theorem}
Let \((X, d)\) and \((Y, \rho)\) be compact metric spaces, \(f_{1,\infty}\) be a sequence of equicontinuous maps from \(X\) into itself, \(g_{1,\infty}\) be a sequence of equicontinuous maps from \(Y\) into itself. If \(\{ \pi_i \}_{i=1}^{\infty}\) is a semiconjugacy from $(X,f_{1,\infty})$ to $(Y,g_{1,\infty})$ and \(\pi_i \in \{ \phi_1, \phi_2, \ldots, \phi_k \},\, k < \infty,\) for every \(i \geq 1\), then for any nonempty subset \(Z \subseteq X\) and any \(\theta \in [0,1]\),
\[
\underline{h}_{\mathrm{top}}\left(f_{1,\infty},Z,\theta\right)
\;\le\;
\underline{h}_{\mathrm{top}}\left(g_{1,\infty},\pi_{1}(Z),\theta\right)
\;+\;
\max_i \sup_{y \in Y} H\bigl(f_{1,\infty}; \pi_i^{-1}(y)\bigr),
\]
\[
\overline{h}_{\mathrm{top}}\left(f_{1,\infty},Z,\theta \right)
\;\le\;
\overline{h}_{\mathrm{top}}\left(g_{1,\infty},\pi_{1}(Z),\theta \right)
\;+\;
\max_i \sup_{y \in Y} H\bigl(f_{1,\infty}; \pi_i^{-1}(y)\bigr).
\]
\end{theorem}
\begin{proof}
Building upon Bowen’s original ideas~\cite{bowen1971entropy} 
and the subsequent extensions by Kolyada and Snoha~\cite{kolyada1996topological}, 
as well as by Fang et al.~\cite{fang2012dimensions}, 
we now establish an extension of these results for $\theta$-intermediate topological entropies.

Let \(a=\max_i \sup_{y \in Y} H\bigl(f_{1,\infty}; \pi_i^{-1}(y)\bigr)\). If \(a= \infty\) there is nothing to prove. So we can assume that \(a< \infty\).

Fix \(\epsilon > 0\) and \(\tau > 0\). For each \( y \in Y \) and \(j \in \{1,\ldots,k\}\), choose \( m_{j}(y) \in \mathbb{N} \) such that
\begin{align*}
a + \tau &\geq \max_{j \in \{1, \ldots, k\}} H\bigl(f_{1,\infty}; \phi_j^{-1}(y); \epsilon\bigr) + \tau \geq H\bigl(f_{1,\infty}; \phi_j^{-1}(y); \epsilon\bigr) + \tau \\
&\geq \frac{1}{m_j(y)} \log r_{m_j(y)}^*\bigl(f_{1,\infty}; \phi_j^{-1}(y); \epsilon\bigr),
\end{align*}
where 
\[
H\bigl(f_{1,\infty}; \phi_j^{-1}(y); \epsilon\bigr)=\limsup_{n\to\infty}\frac{1}{n}\log r_n^*\bigl(f_{1,\infty};\phi_j^{-1}(y);\epsilon\bigr).
\]
Let \( E_y^*(j) \) be a \((m_j(y), \epsilon)^*\)-spanning set of \(\phi_j^{-1}(y)\) with respect to \(f_{1,\infty}\), satisfying 
\[
\mathrm{card}\bigl(E_y^*(j)\bigr) = r_{m_j(y)}^*\left(f_{1,\infty}; \phi_j^{-1}(y); \epsilon\right)
\]
for each \( j \in \{1, \ldots, k\} \).
Define
\[
U_y(j) = \left\{ u \in X : \exists z \in E_y^*(j) \text{ such that } d_{m_{j}(y)}^{*}(u, z) < 2\epsilon \right\}.
\]
Then \( U_y(j) \) is an open neighborhood of \(\phi_j^{-1}(y)\) and 
\[
\left(X \setminus  U_y(j)\right) \cap \bigcap_{\gamma > 0} \phi_j^{-1}(\overline{B_\gamma(y)}) = \emptyset,
\]
where \( B_\gamma(y) = \{ y' \in Y : \rho(y', y) < \gamma \} \). By the finite intersection property of compact sets, there exists \( W_y(j) = B_{\gamma_j(y)}(y) \) for which \( U_y(j) \supseteq \phi_{j}^{-1}(W_y(j)) \). Denote \(W_y=\bigcap_{j\in \{1,\ldots,k\}}W_y(j)\), then \(\phi_j^{-1}(W_y)  \subseteq U_y(j)\) for every \(j \in \{1,\ldots,k\}\).

Since \( Y \) is compact, there exist \( y_1, \dots, y_p \) such that \( \{W_{y_1}, \dots, W_{y_p}\} \) cover \( Y \). Let \( \delta_1 > 0 \) be a Lebesgue number of the open cover \( \{W_{y_1}, \dots, W_{y_p}\} \) with respect to \(\rho\), and set \(0< \delta< \delta_1 / 2 \), \( M:=\max_{1\le j\le k, 1\le t\le p} m_j(y_t)\).

For every \(i \ge 1\) there exists some \(j_{i} \in \{1,\ldots,k\}\) with \(\pi_{i}=\phi_{j_i}\).
Now, for \( y \in Y \) and \( m \in \mathbb{N} \), we claim there exist \( \ell(y) > 0 \) and \( v_1(y), \dots, v_{\ell(y)}(y) \in X \) such that

\[
\ell(y) \leq e^{(a + \tau)(m + M)} \quad \text{and} \quad \bigcup_{i=1}^{\ell(y)} B_{m,f_{1,\infty}}(v_i(y), 4\epsilon) \supseteq \pi_{1}^{-1}(B_{m,g_{1,\infty}}(y, \delta)),
\]
where 
\[
B_{m,g_{1,\infty}}(y, \delta)= \left\{ y' \in Y : \rho_m(y, y') < \delta \right\}
\]
and 
\[
B_{m,f_{1,\infty}}(v_i(y), 4\epsilon)= \left\{ x' \in X : d_m(v_i(y), x') < 4\epsilon \right\}.
\]

Take \(y \in Y\) and let \(c_0(y) \in \{ y_1, \ldots, y_p \}\) be such that \(W_{c_0(y)} \supseteq \overline{B}_{\delta}(y)\) and define \(t_0(y) = 0\). Next let \(t_1(y) = m_{j_1}(c_0(y))\) and let \(c_1(y) \in \{ y_1, \ldots, y_p \}\) be such that \(W_{c_1(y)} \supseteq \overline{B}_{\delta}(g_1^{t_1(y)}(y))\). Similarly, assuming that \(t_0(y), \ldots, t_s(y)\) and \(c_0(y), \ldots, c_s(y)\) are already defined, we define \(t_{s+1}(y) = t_s(y) + m_{j_{t_s + 1}}(c_s(y))\) and let \(c_{s+1}(y) \in \{ y_1, \ldots, y_p \}\) be such that \(W_{c_{s+1}(y)} \supseteq \overline{B}_{\delta}\left(g_1^{t_{s+1}(y)}(y)\right)\). Finally, let \(q =q(y)\) be such that 
\[ \sum_{s=0}^{q-1} m_{j_{t_s + 1}}(c_s(y)) = t_q(y) < m \leq t_q(y) + m_{j_{t_q + 1}}(c_q(y)). 
\]

For \( z_0 \in  E_{c_0(y)}^*(j_1), \dots, z_{l} \in E_{c_q(y)}^*(j_{t_q + 1}) \), define
\[
\begin{aligned}
V(y; z_0, \dots, z_{q}) 
= \Bigl\{\, 
  &u \in X : 
  d\left(f^{t + t_s(y)}_1(u), f_{t_s(y)+1}^t(z_s)\right) < 2\epsilon \\[0.4em]
  &\text{for all } 0 \le t < m_{j_{t_s(y) + 1}}(c_s(y)),~ 0 \le s \le q(y) 
  \Bigr\}.
\end{aligned}
\]

It is not hard to verify that
\begin{equation}\label{eq41}
\bigcup_{z_0 \in  E_{c_0(y)}^*(j_1), \dots, z_{q} \in E_{c_q(y)}^*(j_{t_q + 1})} V\left(y; z_0, \dots, z_{q}\right) \supseteq \pi_{1}^{-1}\left(B_{m,g_{1,\infty}}(y, \delta)\right).  
\end{equation}

For each such tuple \( (z_0, \dots, z_{q}) \), pick any \( v(z_0, \dots, z_{q}) \in V(y; z_0, \dots, z_{q}) \). Then 
\[
B_{m,f_{1,\infty}}\left(v(z_0, \dots, z_{q}), 4\epsilon \right) \supseteq V(y; z_0, \dots, z_{q}).
\]
Moreover, by (\ref{eq41}), we have 
\begin{equation}\label{eq42}
\bigcup_{z_0 \in  E_{c_0(y)}^*(j_1), \dots, z_{q} \in E_{c_q(y)}^*(j_{t_q + 1})} B_{m,f_{1,\infty}}\left(v(z_0, \dots, z_{q}), 4\epsilon \right) \supseteq \pi_{1}^{-1}\left(B_{m,g_{1,\infty}}(y, \delta)\right).
\end{equation}

Let
\[
\ell(y) = \prod_{s=0}^{q(y)} \mathrm{card}\left(E_{c_s(y)}^*(j_{t_s(y) + 1})\right)=\prod_{s=0}^{q(y)} r^{*}_{m_{j_{t_s(y) + 1}}(c_s(y))}\left(f_{1,\infty},\pi^{-1}_{t_s(y)+1}(y), \epsilon\right).
\]

Clearly,
\begin{align*}
\ell(y) &=\exp \left(\sum_{s=0}^{q(y)}\log r^{*}_{m_{j_{t_s(y) + 1}}(c_s(y))}\left(f_{1,\infty},\pi^{-1}_{t_s(y)+1}(y), \epsilon\right)\right)\\
& \leq \exp \biggl((a+\tau)\sum_{s=0}^{q(y)}m_{j_{t_s(y) + 1}}(c_s(y))\biggr)\\
& \leq \exp\left((a+\tau)(m+M)\right).  
\end{align*}

Since the number of permissible \((z_0, \ldots, z_{q(y)})\) is \(\ell(y)\), we may let \(v_1(y), v_2(y), \ldots, v_{\ell(y)}(y)\) be an enumeration of
\[
\left\{v(z_0, \ldots, z_{q(y)}) : z_0 \in E_{c_0(y)}^*(j_1), \ldots, z_{q(y)} \in E_{c_q(y)}^*(j_{t_q + 1})\right\}.
\]
Then by (\ref{eq42}),
\[
\bigcup_{i=1}^{\ell(y)} B_{m,f_{1,\infty}}\left(v_i(y), 4\epsilon \right)  \supseteq \pi_{1}^{-1}\left(B_{m,g_{1,\infty}}(y, \delta)\right).
\]
This proves the claim.

Let \(\left\{B_{n_j,g_{1,\infty}}(w_j, \delta) \right\}_{j=1}^\infty\) be a cover of \( \pi_{1}(Z) \) and \( n \leq n_j < n/\theta+1 \) for each \( j \). By the claim, for each \( j \), there exist \( \ell(w_j) \leq e^{(a + \tau)(n_j + M)} \) and \( v_1(w_j), \dots, v_{\ell(w_j)}(w_j) \) such that
\[
\bigcup_{i=1}^{\ell(w_j)} B_{n_j,f_{1,\infty}}\left(v_i(w_j), 4\epsilon \right) \supseteq \pi_{1}^{-1}\left(B_{n_j,g_{1,\infty}}(w_j, \delta)\right).
\]

Thus,
\[
\bigcup_{j=1}^\infty \bigcup_{i=1}^{\ell(w_j)} B_{n_j,f_{1,\infty}}\left(v_i(w_j), 4\epsilon \right) \supseteq \bigcup_{j=1}^\infty \pi_{1}^{-1}\left(B_{n_j,g_{1,\infty}}(w_j, \delta)\right) \supseteq \pi_{1}^{-1}(\pi_{1}(Z)) \supseteq Z.
\]
Then
\[
M_{f_{1,\infty}}\left(Z,s,4\epsilon,n,\theta \right) \leq 
\sum_{j=1}^\infty \sum_{i=1}^{\ell(w_j)} e^{-s n_j} = \sum_{j=1}^\infty \ell(w_j) \, e^{-s n_j} \leq  e^{(a + \tau) M} \sum_{j=1}^\infty e^{- (s - (a + \tau)) n_j}.
\]
Since the above inequality is true for any \(\left\{ B_{n_j,g_{1,\infty}}(w_j, \delta) \right\}_{j=1}^{\infty}\), we have
\[
M_{f_{1,\infty}}\left(Z,s,4\epsilon,n,\theta \right) \leq M_{g_{1,\infty}}\left(\pi_{1}(Z),s-(a+\tau),\delta,n,\theta \right)e^{(a+\tau)M}.
\]
Taking the $\limsup$ as $n \to \infty$ yields:
\[
\overline{m}_{f_{1,\infty}}\left(Z,s,4\epsilon,\theta \right) \leq \overline{m}_{g_{1,\infty}}\left(\pi_{1}(Z),s-(a+\tau),\delta,\theta \right)e^{(a+\tau)M}.
\]
This implies
\[
\overline{h}_{\rm top}\left(f_{1,\infty}, Z, 4\epsilon, \theta \right) \leq \overline{h}_{\rm top}\left(g_{1,\infty}, \pi_{1}(Z), \delta,  \theta \right) + a + \tau \leq \overline{h}_{\rm top}\left(g_{1,\infty}, \pi_{1}(Z), \theta \right) + a + \tau.
\]
Letting \(\varepsilon\to0\) and \(\tau\to0\) completes the proof for the upper \(\theta\)-topological entropy.  
The argument for \(\underline{h}_{\mathrm{top}}\) is analogous, replacing each \(\limsup\) by \(\liminf\).
\end{proof}

\section{Example}\label{ex}

Let $\Sigma_2=\{0,1\}^{\mathbb Z}$ be the sequence space of all bi-infinite sequences of two symbols, where $\mathbb Z$ denotes the set of all integers.
Let
\[
\sigma_1: \Sigma_2 \to \Sigma_2, \quad
(w_i)_{i=-\infty}^{+\infty} \mapsto (w_{i+1})_{i=-\infty}^{+\infty}
\]
be the classical (one-step) shift and
\[
\sigma_k: \Sigma_2 \to \Sigma_2, \quad 
(w_i)_{i=-\infty}^{+\infty} \mapsto (w_{i+k})_{i=-\infty}^{+\infty}
\]
be a \( k \)-step shift for \( k \in \mathbb{N} \).

The following example is adapted from Pesin~\cite[Example~11.2, p.~76]{pesin1997dimension}.
\begin{example}
Let \( f_1 = \sigma_2 \) and \(f_i = \sigma_1\) for \( i > 1 \). Consider the NDS \(\left(\Sigma_2, f_{1,\infty}\right)\), then we have:
\begin{itemize}
\item[(1)] 
For any \(\theta \in [0,1]\),
\[
\underline{h}_{\mathrm{top}}\left(f_{1,\infty},\Sigma_2,\theta \right)
= \overline{h}_{\mathrm{top}}\left(f_{1,\infty},\Sigma_2,\theta \right)= \log 2.
\]

\item[(2)]  Define
\[
Z_k = \left\{ \omega = (\omega_n) \in \Sigma_2 : \omega_n = 1 \text{ for all } |n| \ge k \right\}, 
\qquad 
Z = \bigcup_{k \in \mathbb{N} \cup \{0\}} Z_k.
\]
Then
$
h_{\mathrm{top}}^{B}\left(f_{1,\infty}, Z\right) = 0,
$
while for any \(\theta \in (0,1]\),
\[
\underline{h}_{\mathrm{top}}\left(f_{1,\infty}, Z, \theta\right)
= \overline{h}_{\mathrm{top}}\left(f_{1,\infty}, Z, \theta\right)
= \log 2.
\]
\end{itemize}
\end{example}

\begin{proof}
(1) Since \(\sigma_1(\Sigma_2)=\Sigma_2\), it follows from Corollary~\ref{momo-theta} that
\[
\underline{h}_{\mathrm{top}}\left(f_{1,\infty},\Sigma_2,\theta\right)=\underline{h}_{\mathrm{top}}\left(\sigma_1,\Sigma_2,\theta\right) \qquad
\overline{h}_{\mathrm{top}}\left(f_{1,\infty},\Sigma_2,\theta\right)=\overline{h}_{\mathrm{top}}\left(\sigma_1,\Sigma_2,\theta\right).
\]
Moreover, since \(\Sigma_2\) is compact and invariant, Theorem~11.5 in~\cite{pesin1997dimension} gives
\[
h^{B}_{\mathrm{top}}\left(\sigma_1,\Sigma_2\right)
= h\left(\sigma_1,\Sigma_2\right)
= \log 2,
\]
which implies the desired result.

(2) It is straightforward that \(\overline{Z} = \Sigma_2\). 
Hence, by Proposition~\ref{sc} together with (1), we obtain for any \(\theta \in (0,1]\),
\[
\underline{h}_{\mathrm{top}}\left(f_{1,\infty}, Z, \theta \right)
= \underline{h}_{\mathrm{top}}\left(f_{1,\infty}, \Sigma_2, \theta \right)
= \log 2,
\]
\[
\overline{h}_{\mathrm{top}}\left(f_{1,\infty}, Z, \theta \right)
= \overline{h}_{\mathrm{top}}\left(f_{1,\infty}, \Sigma_2, \theta \right)
= \log 2.
\]
On the other hand, Example~11.2 in~\cite{pesin1997dimension} shows that
\(
h_{\mathrm{top}}^{B}(\sigma_1,Z)=0.
\)
Since $\sigma_1(Z)=Z$, Corollary~\ref{momo-theta} implies that
\[
h_{\mathrm{top}}^{B}\bigl(f_{1,\infty},Z\bigr)
=
h_{\mathrm{top}}^{B}(\sigma_1,Z)
=
0.
\]
Therefore, the functions 
$\theta \mapsto \underline{h}_{\mathrm{top}}\left(f_{1,\infty}, Z, \theta\right)$ 
and 
$\theta \mapsto \overline{h}_{\mathrm{top}}\left(f_{1,\infty}, Z, \theta\right)$ 
are both discontinuous at $\theta = 0$ but constant for all $\theta \in (0,1]$.
\end{proof}

\section*{Acknowledgements}
The author is sincerely grateful to the anonymous referee for the valuable comments and suggestions that greatly improved the quality of this manuscript.

\section*{Declarations}

\section*{Funding}
This work was supported by the Science and Technology 
Research Program of Chongqing Municipal Education Commission (Grant No.KJQN202500802).

\section*{Competing Interests}

The author declares that there are no conflicts of interest regarding this paper.

\bibliography{references}

\end{document}